# The Ethical and Sustainable Concerns Triangle:
# A Framework for Navigating Discourses in Mathematics and Its Education

**Dennis Müller[1]**     **Maurice Chiodo[2]**

**Michael Meyer[3]**

**13 October 2025**

**Abstract:** The literature on ethics and sustainability in mathematics and its education is increasingly complex and fragmented, potentially leading to communication breakdowns between different scholarly traditions. To address this, the paper introduces the "Ethical and Sustainable Concerns Triangle," a framework that maps discourses based on their relative concern for three areas (represented as three vertices in the triangle): "Mathematics", "Community", and "Society/Planet". By integrating a systems theoretic perspective, we analyse discourses as dynamic systemic reactions to external irritations. Our analysis reveals that the field's fragmentation can be explained by the "location effect": the phenomenon whereby a discourse's position within the triangle shapes its perception and acceptance of other scholarship. By mapping key discourses and educator archetypes, the framework functions as a meta-heuristic tool. Ultimately, it serves not only to facilitate critical reflection but also as a call for the epistemic humility and dialogue needed to advance the field.

---

[1] Institute of Mathematics Education, University of Cologne, Germany. dennis.mueller@uni-koeln.de

[2] Centre for the Study of Existential Risk, University of Cambridge, U.K. mcc56@cam.ac.uk

[3] Institute of Mathematics Education, University of Cologne, Germany. michael.meyer@uni-koeln.de

# Table of Contents

# Introduction

Six years ago, some of the authors organised the *2nd Meeting on Ethics in Mathematics*.[4] The conference brought together scholars from different fields across the globe (mathematics, education, philosophy, psychology, etc.) who were all concerned with aspects of ethics, social justice, and sustainability in mathematics and its education. However, the anticipated productive discussions largely failed to materialise; instead, the conference fractured into competing camps. Rittberg, who presented a case study on epistemic injustice, later recalled the events:

> "[S]ome mathematicians got quite upset […] It was a conference that had a schism and was by far the most exciting conference I ever went to." (Rittberg, as cited in Müller, 2024, p. 72)

What happened at the conference was a clash between traditionalist positions (on philosophy, mathematics, ethics, and education) with critical positions (critical mathematics education, mathematics for social justice, and critical pedagogy at large) (Müller, 2024). Unspoken assumptions about the nature of mathematics, education, and ethical concerns came to light in a schismatic way and caught many of the attendees by surprise. This event signalled to some of the authors how rapidly the landscape of ethics in mathematics and its education descends into complexity, particularly when diverse international stakeholders are involved. This rise in complexity is related to a recent surge in publications as "mathematics education has arrived late to surfacing questions around sustainability" (Makramalla et al., 2025) and ethics (Müller, 2025a), and due to the socio-cultural, economic, and political complexity accompanying these questions (e.g., Amico et al., 2023) — whose ethics (Boylan, 2016) and sustainability (Makramalla et al., 2025) are we considering?

This paper continues the initial analysis presented by Müller (2024) and posits that to advance the conversation, a new heuristic analytical perspective which can further map and explain some of the fragmentation would be beneficial.[5] We therefore introduce the "Ethical and Sustainable Concerns Triangle" (the ESCT) as a framework to analyse key dimensions of ethics and sustainability in mathematics and its education. The development of this framework aims to synthesise the knowledge gained from various systematic and non-systematic reviews, as well as from some of the area's key papers, in a structured and visual way. To develop these perspectives, we adjust the areas' three prototypical ethical concerns (Chiodo & Müller, 2024, p. 940) to include sustainability:

1) the ethics and sustainability concerns of ensuring that (research) mathematics remains a stable and continuing body of knowledge,
2) the ethical, sustainable, and social questions surrounding our mathematical, and other, related communities, and
3) questions regarding mathematics within, and its impact on, wider society and our planet.

---

[4] https://www.ethics-in-mathematics.com/EiM2/

[5] Enabling communication between different discourses also ensures that methods, tools, and analyses are not accidentally and independently reproduced and/or discovered.

We construct a triangle whereby the concerns (1) - (3) become three vertices, which we label "Mathematics", "Community", and "Society/Planet" (Figure 1). Here, we use "Society/Planet" rather than "Society/Ecological" to also capture discourses that worry about planetary scales (e.g., Burden et al., 2024, 2025; Crawford, 2021), to recognise that sustainability (in mathematics education) is complex and contested regarding different human and non-human interests (cf. Makramalla, 2025), and to capture recent desires for an emphasis on bigger sustainability questions (Bakker et al., 2025). We use "Community" rather than "Human" to show that the concerns about the mathematical community can be broader than educational concerns or concerns about the people involved. They can also ask about meta-developments in the mathematical community, by asking who a mathematician is (cf. Buckmire et al., 2023; Müller et al., 2022) or by being concerned about dynamics of epistemic injustice within mathematics education (e.g., Tanswell & Rittberg, 2020; Rittberg et al., 2020; Rittberg, 2020, 2023, 2024; Hunsicker & Rittberg, 2022).

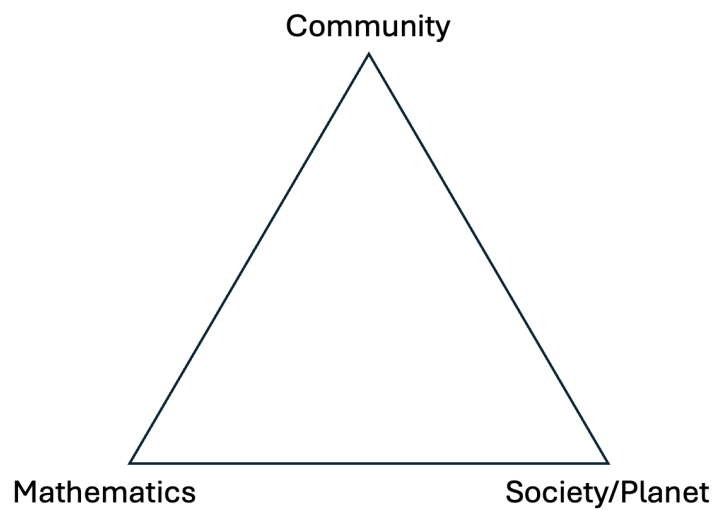

Figure 1: The Ethical and Sustainable Concerns Triangle (ESCT)

Onto this (equilateral) triangle, we map different aspects of the literature, whereby the proximity to a vertex is understood as a degree of concern, motivation or focus. For example, a discourse that is close to the "Society/Planet" vertex may largely find its motivations in or focus on larger societal or planetary concerns. The idea behind this is that convex analysis[6] allows us to represent each point in the triangle as a convex combination of the vertices (Boyd & Vandenberghe, 2004).

**Definition:** A convex combination of a triangle's three vertices $x_1, x_2, x_3$ is given by $\lambda_1 x_1 + \lambda_2 x_2 + \lambda_3 x_3$ where $\lambda_1, \lambda_2, \lambda_3 \geq 0$ and $\lambda_1 + \lambda_2 + \lambda_3 = 1$.

**Theorem:** All points of a triangle can be represented by a convex combination of the three vertices.

In this sense, the vertices represent idealised concerns.[7] In reality, most people will have concerns which can be understood as covering multiple vertices, and thus are understood through (weighted) combinations in this framework, rather than as actually "touching" different vertices - a different metaphorical picture. In our setting, the convex combinations represent how much of the focus lies on a specific vertex. A higher weight for a specific vertex represents a position that gives more importance to concerns associated with this vertex.[8] Consider, for example, an approach to inclusive classrooms ("Community" vertex, "$\lambda_{community}$ is large") that is also a response to societal injustices ("Society/Planet" vertex, "$\lambda_{society}$ is large"); the triangle would illustrate this as a discourse located in the space between those two vertices.

Our analysis will demonstrate that the apparent fragmentation of the field can be understood as a series of relational positions and dynamic tensions between these three core concerns, represented by the three vertices of the triangle. Convex analysis gives us the inspiration to consider convex combinations of concerns. However, given the interpretative acts involved in any (discourse) analysis and mapping of concerns, we suggest to the reader that it can sometimes be better to think of regions where specific concerns are strong, rather than to think that the exact mathematical precision from convex analysis transfers over to the ESCT (as illustrated in Figure 2).

---

[6] In the mathematical sense.

[7] Here we are using three concerns, but this would also work for more concerns by considering an n-dimensional simplex. The methodology developed in this paper can then also be applied in this new context.

[8] The literature on fairness (e.g., Pfeiffer et al., 2023; Mitchell et al., 2021) has mathematically proven that different (formal) definitions of fairness and justice can be mutually incompatible. Therefore, trade-offs are a necessary reality, and not an option - reflected here through the convex combinations of different concerns.

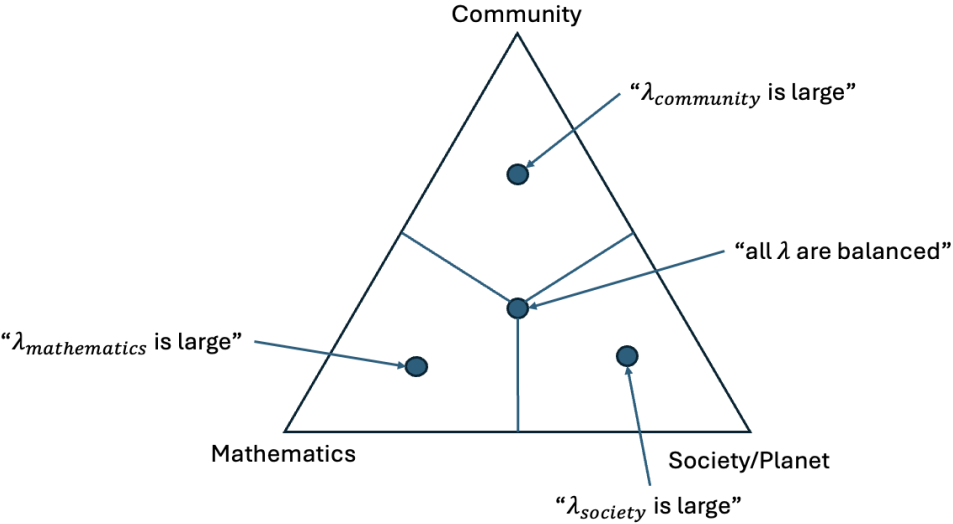

Figure 2: Convex combinations of concerns

So how does the ESCT present itself in practice? This is perhaps best illustrated with a concrete example. In early 2025, one of the authors gave a talk to an audience of academic mathematicians on teaching ethics and embedding it in everyday (university) mathematics teaching - a talk firmly rooted in the "Society/Planet" vertex of the ESCT. During the Q&A session after the talk, one mathematician asked (paraphrased), "How can I fit this into my teaching when the curriculum is already packed with mathematics? I have no spare space or time for this." - a concern firmly rooted in the "Mathematics" vertex of the ESCT. Immediately afterwards, another mathematician asked (paraphrased) "I'm concerned that this teaching might affect how my students seek out jobs, and that on ethical grounds they might choose not to apply for certain types of work, thus reducing their opportunities." - a concern firmly rooted in the "Community" vertex of the ESCT. To the first question, the speaker replied with "Rather than optimising only over the (technical) mathematical content being taught, might you consider reducing that *slightly* to incorporate some ethics teaching; material students might otherwise never see at all?" To the second question, he replied with "Rather than optimising only over the interests and careers of your students, might you consider adjusting that *slightly* to incorporate the interests of those who your students may impact in society through their mathematical work?" Both answers were well-received by those posing the questions, and the wider audience. Here, the ESCT helped to quickly identify different concerns, acknowledge their individual validity, and then seek out the possibility of a compromise.

The ESCT builds on a long tradition of using triangles as explanatory devices in mathematics education. These models, however, typically focus on the pedagogical situation itself. The didactic triangle, for instance, describes the relationships between teacher, student, and content (Goodchild & Sriraman, 2012), and epistemological triangles map the cognitive construction of knowledge (Steinbring, 1998). The model of didactic reconstruction serves as a framework for designing learning environments by systematically connecting three core components: subject-matter clarification (fachliche Klärung), recording of learners' perspectives (Erfassung von Lernerperspektiven), and didactic structuring (didaktische Strukturierung) (e.g., Kattmann, 2007; Kattman et al., 1997; Reinfried et al., 2009), including in mathematics (Prediger, 2005). The ESCT departs from this direct focus on teaching practice. It operates on a meta-level as a heuristic tool designed to analyse scholarly discourses within the field. Its vertices represent idealised concerns that orient researchers and educators, and thus also pedagogical actors or design processes. The ESCT is therefore a tool for critical reflection on the field itself, not a direct guide for designing or evaluating individual learning situations.

In proposing the ESCT, the authors are mindful of the critiques levelled against the didactic triangle. Gerwig (2017) warned of a "dangerous [didactic and analytical] shortening" when he argued that the elegance of the didactic triangle risks oversimplifications and begins to obscure the very classroom situations it wants to enlighten. Taking on lessons from responsible mathematical model building (Chiodo & Müller, 2025; Thompson, 2022), the proposed framework attempts to mitigate this by mapping different aspects onto the same triangle rather than presenting a single visual, which we outline in detail throughout this paper. Given the profound complexity of the professional and academic spheres of mathematics and its education, no single 2D "model" can capture their entirety. The ESCT should therefore be seen as an analytical lens, not a comprehensive map. It allows for various "projections" that illuminate the distinct ways different groups of practitioners approach issues of sustainability and ethics. Each projection necessarily highlights certain details and loses others, but taken together, they provide what we see as a useful meta-heuristic way of navigating the landscape. In doing so, the framework complements Vygotsky's triangle on micro-perspective on mediated actions (Zittoun et al., 2007) by giving a macro-level perspective on mediations between systemic structures of the field itself.

Any model is a deliberate simplification, and the ESCT is no exception. In taking on Chiodo and Müller's (2024) three concerns, we propose that path-dependencies, power dynamics, and the socio-historical context are better understood as critical lenses that can be applied to the map, rather than as additional vertices. Applying additional concepts to the map reveals forces that shape the landscape and causes discourses to cluster in particular regions.

The ESCT is neither designed to be a pedagogical tool, nor is it designed to be a normative map. The intention is to describe, analyse, and enable open forms of communication. By placing various discourses on the same triangle, we do not intend to create an equivalence between different positions or discourses. A position on the map is an interpretation of central concerns and does not imply an endorsement or criticism. But even though it is not explicitly designed to be a pedagogical tool, the framework still carries practical implications beyond academic discourse analysis. For practitioners, the framework can be a tool for (critical) self-reflection. A teacher's awareness of their own location in the triangle may influence their pedagogical choices and their ability to engage with diverse students and complex curricular demands (cf. Müller, 2025a, 2025b; Rycroft-Smith et al., 2024; Chiodo & Bursill-Hall, 2018, for different discussions of ethical awareness).[9]

To answer our research question of developing the ESCT as a heuristic tool to manage discursive complexity in ethics and sustainability in mathematics and its education, we proceed in two parts.

In part one, we discuss some of the foundational aspects of our proposed triangle. We first focus on interpretative acts involved in setting up the triangle as a framework to understand some of the decisions we have made in its setup and their consequences. We will then consider a systems-theoretic interpretation of its vertices to introduce a perspective which will later allow us to speak about the origins of certain discourses. After this, we introduce the "location effect", i.e., how the location inside the triangle and the distance between two locations can affect what a particular discourse (or an individual scholar) can see and accept readily. We will specify this by considering different levels of learning, represented by education about, for, and as sustainability or social justice and specific archetypal educators. In this context, we will also discuss the differences between ethics and sustainability in school and higher mathematics. Having established these foundations, we map the "correlation" between the vertices and neutrality, purity, and universality assumptions.

These foundational perspectives provide the theoretical underpinnings of the second, and final, part of the paper, where we map different ethical and sustainable discourses onto the triangle according to (what we perceived as) their dominant concerns. We will combine these with a discourse analytic lens inspired by Laclau & Mouffe (2020) to establish a methodology to perform this mapping. Different discourses will be represented by discourse clouds, highlighting their spread of concerns.

---

[9] Guiding questions are presented in Appendix C (Table 7). Appendix B also provides a one-page summary of the framework.

# Part I: Foundations of the ESCT

## Interpretative Acts and their Critique

To avoid unnecessary (and misleading) mathematical precision, we will now first speak about some of the interpretative decisions involved in setting up the ESCT. Developing our proposed framework involves different such acts, most notably,

- the decision to turn Chiodo & Müller's (2024) concerns into a triangle, and thereby reducing discursive complexity using simple visuals, and
- the choices made when mapping different aspects (e.g., discourses) onto the triangle.

The choice of the map ("triangle") and the choice of the mapping ("location") are interpretative decisions made to explain tensions and similarities. The authors acknowledge this openly and have tried to make biases visible throughout the text.[10] The structure of the paper is a conscious decision: rather than interpreting different vertices after each other, we interpret them simultaneously from different perspectives. We believe that this serves two important goals: (a) to make our thinking process and thus some of our biases more visible, and (b) to make our approach more accessible to an interdisciplinary audience. This approach increasingly adds depth and nuance through a series of (layered) perspectives which come together in the discursive clouds of Figure 15.

Our proposed triangle can be understood as an attempt to find an abstract, locational grammar that brings to light collective conventions, rules and concerns. In this framework, meaning is not derived in isolation but from the relational position within the triangle. Chiodo's and Müller's (2024) conception of three areas of concern can be critiqued as providing a static snapshot, fixing "Mathematics," "Community," and "Society/Planetary" as universal signifiers rather than treating them as evolving socially constructed discursive concepts. The triangle necessarily also provides a static, interpretive snapshot, and at times reveals this very fragility when it maps different dynamic discourses about mathematics. Ernest notes that

---

[10] Figure 7 and Appendix A go deeper into some of the authors' positions.

"[s]ome authors like to distinguish between what is termed Mathematics and mathematics (Bishop, 1988). Mathematics with a capital M is the sum or body of formal mathematical knowledge and mathematical objects that exist independently of mathematicians. The lower case m in the term mathematics signifies that it refers to and comprises human mathematical practices. The distinction is an ontological one, about existence. Mathematics with a capital M is claimed to be real, made up of independently existing concepts, objects, proofs, truths, and theories." (Ernest, 2020b, p. 5)

Regarding this distinction, the triangle's "Mathematics" vertex errs on the side of Mathematics rather than mathematics. This decision is a deliberate methodological choice, designed to surface the inherent tensions between traditionalist mathematicians' perspectives and their philosophies of mathematics (e.g., Platonism, formalism, and logicism) and more recent developments in 20th and 21st-century philosophy of mathematics and mathematics education (e.g., social constructivist philosophies, critical mathematics education). It enables us to make visible the tensions, antagonisms, and dynamics that different discourses exhibit. This choice shapes the entire framework, leading to competing interpretations for the location of ethnomathematical discourses (cf. Figures 15 and 16). The authors acknowledge that this framing may inadvertently reinforce the very tensions they try to explain, and thus appeal to the reader to be mindful in the usage of the framework. The dynamics and human practices that constitute mathematics come back into play through the ongoing system-theoretic interpretation that happens throughout the paper. It will help to overcome some of the inelasticity of a purely triangular perspective, and help to see both Mathematics and mathematics.[11] Additionally, our perspective regarding "pure" and "applied

---

[11] A radical post-structuralist might see this as evidence to reject the notion of "Mathematics" as a stable vertex altogether, rendering the map (almost) meaningless from their perspective. The authors invite such critiques, but acknowledge that a picture with a "floating 'Mathematics' vertex" could hinder the stated goal of providing a comparatively simple map that fosters dialogue between different stakeholders in mathematics and its education. The decision to pick Laclau and Mouffe (2020) to build the later discourse and analytic mapping is an explicit acknowledgement to post-structuralist thought, and the understanding that mathematics and its education is a complex web of discourses, politics, power, identity, and competing hegemonic claims (and thus, in essence, by a wish to balance structuralist and post-structuralist perspectives). It's also an acknowledgment that meaning largely comes from discourse. Epistemologically, the choice of [M]athematics privileges an abstract ontology of mathematics, which leads to locating ethnomathematics (which is deeply concerned with the nature of mathematical knowledge itself) as being distant from the "Mathematics" vertex. Politically, this placement replicates and reinforces the very marginalisation of culturally diverse mathematical traditions that ethnomathematics is challenging, which motivated the authors to end the paper with an

mathematics" is that this distinction is a question of focus: pure mathematics mainly asks inner-mathematical questions aimed at producing more mathematics, while applied mathematics mainly asks questions aimed at solving or learning about external problems. Thus, we follow the distinction made by Müller (2018), which leads us to having "one Mathematics vertex." Aspects concerning "pure" vs "applied" are relegated to the mapping of different aspects onto the triangle, and those with a "pure" focus will often be closer to the "Mathematics" vertex.

## A Systems-Theoretic View of the Triangle

At various points throughout this paper, we use systems theory (Luhmann, 1987) to analyse certain dynamics of positions within the triangle. Systems theory provides an analytic lens allowing us to treat the triangle not as a set of static concerns, but as a representation of dynamic, interacting systems and their shared environment (see Figure 3). To apply systems theory, we re-[12]interpret the three concerns, and thus the three vertices of the triangle. The concerns outlined by Chiodo & Müller (2024) are about knowledge, practice, impact, etc., and thus about the functioning of specific systems, and their impact on each other and their relationship to their (shared) environment. When moving from a concern to the system and/or environment this concern is connected to, the perspective moves to a dynamic perspective on the people, communications, institutions, etc. involved in the formation and articulation of this concern. This transition carries with it certain simplifications (e.g., it hides some of the internal dynamics of (contested) actor-networks within each system), and is done to provide additional insights about the ESCT. It is best understood as an additional analytic lens, rather than as a projection (in the mathematical sense) which maintains (all) necessary properties. Table 1 provides a summary of the key terminology.

More concretely, the "Mathematics" vertex represents a functional subsystem of the larger scientific system, while the "Society/Planet" vertex represents the broader social, political, economic, and planetary (e.g., ecological) environment that continually sends out irritations. The scientific (sub-)system of mathematics reproduces itself by generating (or discovering)

---

alternative placement of ethnomathematics that is closer to understanding the [M]athematics vertex as a [m]athematics vertex. We will come back to this at the very end of the paper when we connect the mapping of ethnomathematics to a triangular prism as an extension to the triangle which may capture further qualitative nuances.

[12] The first act of interpretation was turning different idealised concerns into the vertices of the triangle, and now we are re-interpreting these concerns and vertices as systems/the environment of those systems.

proofs, theorems, standards of rigour, etc. (cf. Watson, 2023; also cf. Greer et al., 2024; Hersh, [R.], 1991, 2013; Ernest, 1991). For example, the mathematician Andrew Wiles has publicly expressed concerns about the applications of mathematics in finance and how it could damage the reputation of mathematics itself, saying that "[o]ne has to be aware now that mathematics can be misused and that we have to protect its good name" and that mathematics research and its applications can move "towards goals that you might not believe in" (Andrew Wiles, quoted by Silvera (2013)).

Following Watson's (2023) initial application of systems theory to mathematics education, the "Community" vertex is conceptualised not as a container for a predefined set of individuals, but rather as the social system of communication that constitutes the field of mathematics education - *here, understood broadly to include school, university, professional education, etc.* As an example, consider Walk's and Tractenberg's (2025) deck of cards to teach students about ethics in mathematics. Their concerns about "help[ing] students move to higher levels of achievement" (ibid., p. 1) regarding ethical reasoning are reflected in a new discursive practice (a deck of cards) to teach mathematics students about ethics. The "Community" system's function is to translate, refine, and transmit mathematical knowledge. Its communicative acts include classroom discourses, curriculum designs, and discourses about the community (e.g., the agency of students). Depending on their communicative practices and concerns, individuals can transition between both systems.

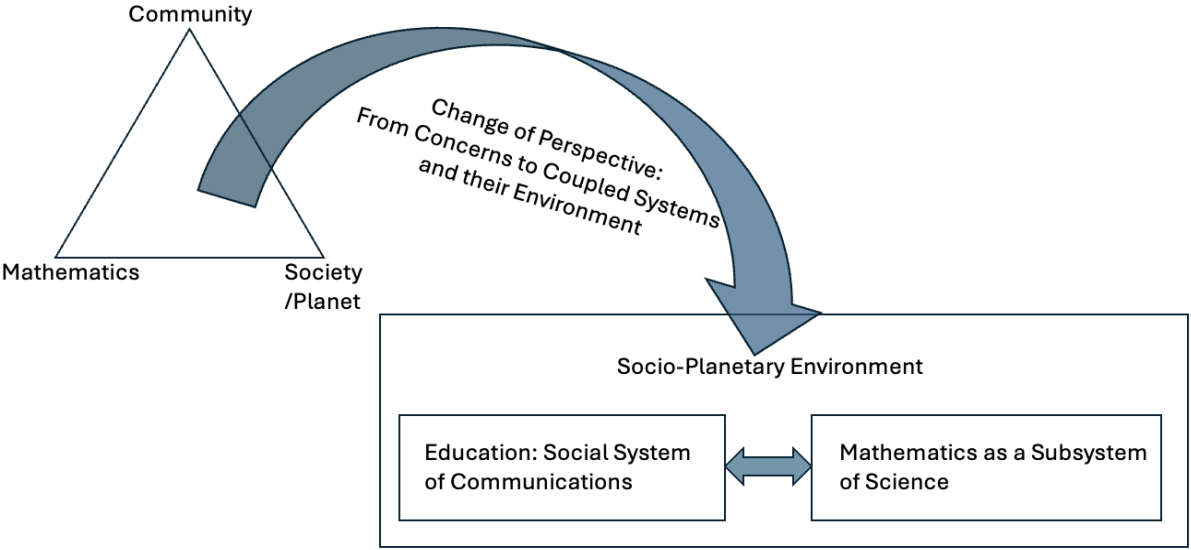

Figure 3: Moving from concerns to systems

As Watson (2023) observes, a system-theoretic perspective of mathematics education carries two profound implications. First, he observes that the education system is operationally closed and driven by autopoiesis - its drive for self-production and self-maintenance. The system is structurally coupled with other systems: it co-evolves through mutual influence without losing its distinct identity. Second, and most crucially, Watson argues that this social system of communication is structurally coupled with the psychic systems (i.e., the consciousness) of the learners, teachers, education scholars, and the many other people involved in its functioning. This explicitly connects the macro-level of social communication to the micro-level of individual cognition, experience, and consciousness.

| Summary of System Theoretic Terminology | | |
|---|---|---|
| **Terminology** | **Definition** | **Usage** |
| Autopoiesis | Self-reproduction of a system when operating according to its internal rules and logic ("self-maintenance"). | Used to describe how mathematics and education maintain their internal consistency in the presence of external pressures. |
| Irritation | Disturbances, demands or pressures that are put onto a system from outside, which cannot be ignored by the system. Irritations are interpreted through the system's internal logic. | Examples: Societal concerns (e.g., climate change, injustices, algorithmic biases) acting as irritations triggering discursive shifts. |
| (Structural) Coupling | Systems are coupled with each other and connected to the environment. Thereby, systems and the environment can influence each other. | This perspective helps us to explain how both the education system and mathematics respond to each other and to the environment. |
| Communication System | Social communication systems are constituted by patterns of communication, rather than by people or institutional structures. | The "Community" vertex is viewed as a social communication system that (re-)produces pedagogical norms, etc. |

Table 1: Summary of system theoretic ideas and their usage

Watson (2023, p. 127) further argues that "mathematics education is, through policy, practice and scholarship in a constant process of defining itself in relation to environmental prompts." External pressures from the socio-planetary environment can irritate the other two systems, and these systems interpret and process irritations according to their own internal logic. As an example, consider the concerns expressed by Alcantara (2024) in her article on California's educational reform plans *Math Education Needs Reform. It Got a War Instead, in* which she notes how mathematicians and computer scientists got worried about new data science courses that do not sufficiently focus on mathematics, while proponents of the reform saw the need in educating underprivileged and minoritised groups.

As we will see, many discourses on ethics and sustainability can be understood as the visible evidence of the educational and mathematical systems processing deep and persistent environmental irritations. Thus, systems theory provides us with a perspective on the origins of a discourse, and allows us to connect the location inside the triangle not just with specific concerns but with the systemic irritation that may have started (or at least strongly influenced) a particular discourse and its approach. A discourse situated between two vertices then represents how one system is irritated by another system or the socio-planetary environment (or how two systems irritate each other). This perspective is necessary to understand why specific discourses may devalue certain concerns, as they are not necessarily perceived as relevant to their fight.

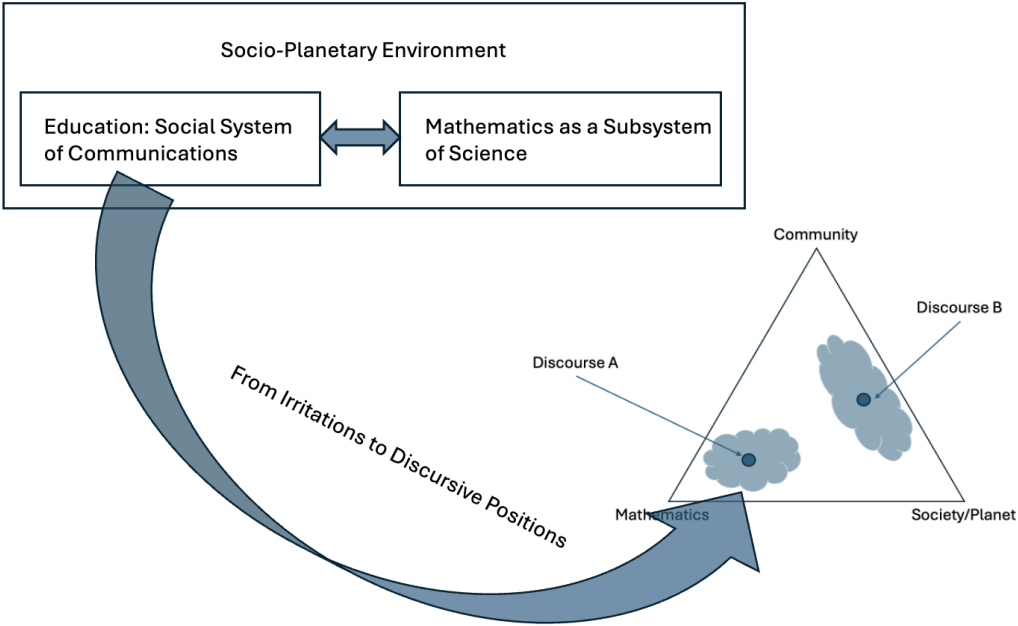

Figure 4: From irritations to discursive positions

For a concrete example, consider societal irritations given by worries about biased algorithms and artificial intelligence, and how they are processed differently by the two systems. "Mathematics" may respond by mathematising fairness and developing fairness-aware algorithms (e.g., Dwork et al., 2012; Kleinberg et al., 2018; Mitchell et al., 2021; Pessach & Shmueli, 2023; Wang et al., 2022; Wong, 2020), while the social system of education may respond by creating new learning modules on data ethics or responding to algorithmically-related equity issues in the classroom (e.g., Baker, [R. S.] et al., 2022; Kasinidou et al., 2021; Jiang & Pardos, 2021). Some individuals respond in both ways (e.g., Porter, 2024).

While likely neither the sociopolitical (Gutierrez, 2013) nor the ethical turn (Müller, 2025a; Ernest, 2024a) of mathematics and its education represent full paradigm shifts according to Kuhn's (1996) *Theory of Scientific Revolutions*, the system-theoretic perspective suggests that a non-negligible number of irritations are happening to traditional perspectives on mathematics. From this perspective, the era of normal mathematics and its education appears to be slowly breaking down, as the shared set of beliefs and theories diminishes. These turns can be understood as early reactions to anomalies which can no longer be suitably addressed by a perspective that positions mathematics as neutral, pure and universal, suggesting that some of the discourses that we will map onto the triangle are turning into competing paradigms in a Kuhnian sense, as also evidenced by the conflict from the 2nd Meeting on Ethics in Mathematics that started this paper or the long-lasting aftermath of the US' math wars (cf. Schoenfeld, 2004; Klein, 2007; Woods, 2023). These wars are not just about the correct ways of *teaching* mathematics, but fundamental debates about what mathematics actually *is*. As such, the ESCT might serve a useful purpose in these debates, which is to help participants see the existence of, and perhaps even some merit to, other paradigms and positions.

## The Location Effect

Scholars continuously make decisions about the literature they engage with, and what they perceive as scientifically valuable, acceptable, and correct. This and the next subsections explore how these decisions are related to the distance between different locations inside the triangle. The location and distance of a position relative to other locations highlights how a position (e.g., their beliefs, values, background, and theoretical commitments) can potentially influence their perception, interpretation, and acceptance of different types of scholarship within ethical and sustainable mathematics and its education (cf. Chiodo & Müller, 2024, forthcoming; Müller, 2024, 2025a; Chiodo & Bursill-Hall, 2018, 2019) – a phenomenon which we will coin the "location effect".

From a systems-theoretic perspective, a psychic system (an individual's consciousness) is structurally coupled to other social systems and communications originating from distant locations may not resonate with the system's internal structures and are thus processed as "noise" rather than meaningful information. What one deems "good", "useful", "rigorous", or "valid" might be questioned, dismissed, or accepted relative to one's location and their distance to a different position. For example, someone firmly rooted in traditionalist schools of mathematics education may struggle to accept the premises, methods, and perspectives of critical mathematics education, which views mathematics as a politically charged and socially constructed practice.[13] The circles and ellipses in Figure 5 represent a visual metaphor of what (someone at) a specific location in the triangle may "see"[14] and "accept"; something outside of the ellipse might be understood as noise. In reality, the field of view is likely not circular or elliptical, but rather depends on different discourses, individual characteristics of researchers, power dynamics, institutional constraints, and other characteristics.

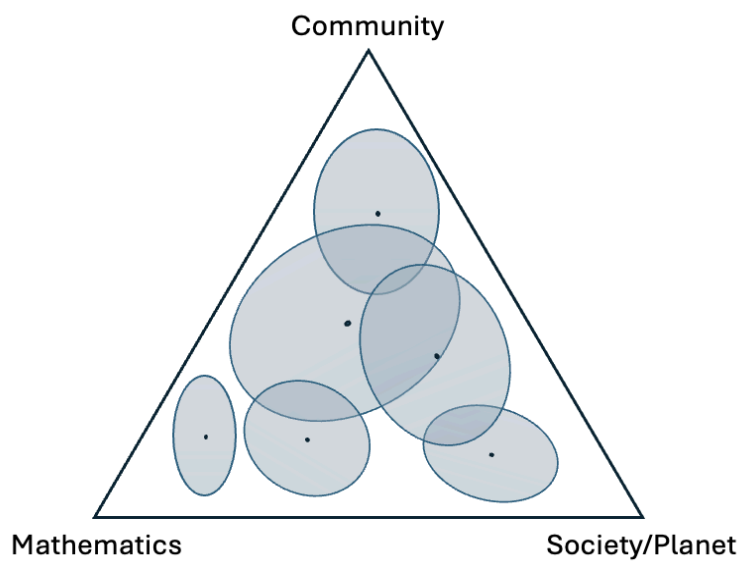

Figure 5: An abstract depiction of the Location Effect[15]

---

[13] As further evidenced both by Crowell's (2022) depiction of school reforms in California, and by Müller's (2024) discussion of the differences between Mathematics for Social Justice and the (Cambridge University) Ethics in Mathematics Project.

[14] Figure 5 is purely metaphorical, without reference to any specific aspect, activity, or event. We give several explicit examples of such diagrams later in the paper, with reference to specific topics.

[15] The authors briefly wondered whether discourses in the middle of the triangle necessarily have a larger field of view. However, they could not conclusively say so, and invite further studies on this.

Studying how mathematics students engage with ethics has shown that concerns for specific aspects of mathematical practice and the vertices of the ESCT are correlated (cf. Chiodo & Müller, forthcoming; building on Chiodo & Müller, 2025), i.e., those who are closely associated with the "Mathematics" vertex typically also deeply care about handling data and information, data manipulation and inference, and the mathematisation of problems, while the "Community" vertex is more often linked to concerns about diversity and perspectives, to the politics of mathematics, and its communication. This can create tensions between the three concerns and the respective systems and their environment (Figure 6), as the needs of one system (e.g., the mathematicians' need for more technical mathematics in teaching) can be seen as an irritating external demand by another system (e.g., an educational community that wants to put more focus on student well-being).

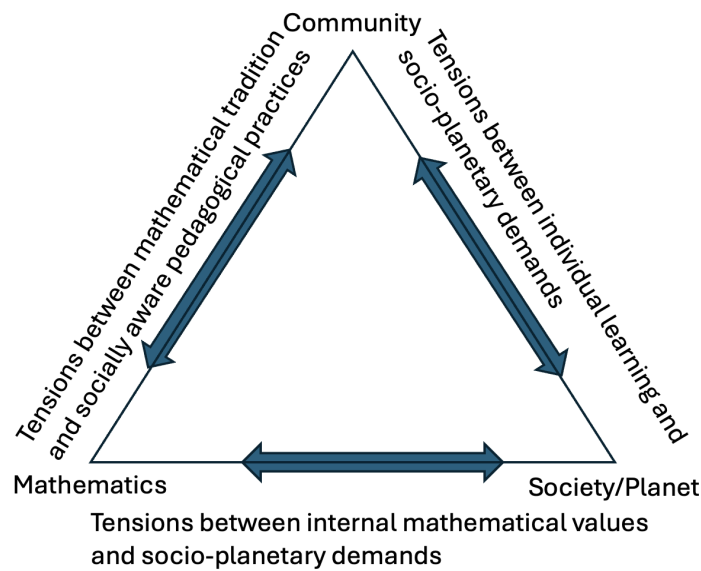

Figure 6: Concern-related and systemic tensions

Behind the location effect are also sociomathematical norms. These norms differ from purely social norms (e.g., the expectation to be a good student or to be a good researcher) in that they are related to activities and interactions within mathematical contexts. The literature has shown that these norms are not fully pre-determined but are continuously (re-)negotiated and regenerated through students', teachers' or researchers' discursive practices and their interactions (e.g., Meyer & Schwarzkopf, 2025; Yackel & Cobb, 1996; Stephan, 2020). The location effect, and in particular, the experience of processing an external claim as "noise", can thus be understood as different groups operating with different sets of sociomathematical norms. The fragmentation of the field and the conflicts between discourses can be understood as processes of re-negotiation.

Epistemic injustices become particularly relevant here. Rigidly adhering to one perspective can lead to the devaluation of legitimate ways of knowing and knowledge claims, effectively functioning as gatekeeping through epistemic injustices and exclusions in mathematics and its education (e.g., Tanswell & Rittberg, 2020; Rittberg et al., 2020; Rittberg, 2020, 2023, 2024; Hunsicker & Rittberg, 2022). Openness, humility, and the bravery to consider alternative claims are central to developing ethics and sustainability further (Chiodo & Müller, 2024; Müller, 2025a).[16] Hence, we want to make explicit that the proposed framework goes beyond merely acknowledging the existence of (locational) biases, and that it also represents a call for epistemic humility, dialogue, and perspective. The complex, multifaceted nature of the Sustainable Development Goals (United Nations, n.d.) requires scholars and educators to adopt nuanced, well-reflective perspectives that look beyond their immediate vicinity.

To further illustrate this and to explain where some of the authors are coming from, Figure 7 shows the location effect for their pragmatic approach to ethics in mathematics. Two of the authors (both trained mathematicians) became interested in the ethics of mathematical practice when they saw the way that mathematics was used and misused in unsustainable and unjust ways in wider society (e.g., Chiodo & Müller, 2018, 2020, 2025; CUEiMS, n.d.; EiMP, 2025). Consequently, these authors initially focused strongly on socio-planetary concerns, attempting to develop pragmatic tools that can help mathematicians to become more ethically aware. The strong focus on socio-planetary concerns meant that, at least initially, there was a strong concern for educating those who know mathematics and are already wielding its power in society, rather than on other community concerns. Hence, as depicted in Figure 7, concerns satisfying both $\lambda_{community} \gg \lambda_{mathematics}$ and $\lambda_{community} \gg \lambda_{society/planet}$ were at various points difficult to process from their location without a corresponding change of (theoretical) perspective, e.g., from pragmatism to critical pragmatism (Müller, 2025b).

---

[16] Fricker (2007) observes that epistemic injustice comes in testimonial and hermeneutical forms. Testimonial injustice occurs when reasonable claims are deemed less credible because of pre-existing prejudices. Hermeneutical injustices refer to situations whereby one party lacks the knowledge and intellectual concepts to adequately communicate their claims to another party (who may have a different knowledge set and understanding of concepts).

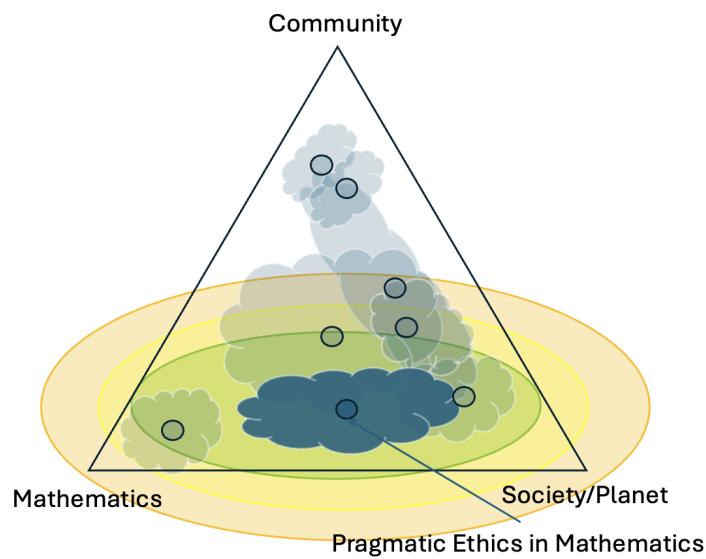

Figure 7: The location effect for pragmatic ethics in mathematics[17]

We will now look at different aspects of the literature on ethics and sustainability in more detail, starting with the distinction between mathematics education about, for, or as sustainability and social justice; different archetypal educators; and differences between school and higher mathematics.

## Mathematics Education *about*, *for* or *as* …?

Connected to the location effect are three distinct approaches to ethical and sustainable education that have emerged, and have thus far remained central for many authors (cf. Makramalla et al., 2025). These can be roughly characterised as

1. mathematics education **about** social justice and sustainability,
2. mathematics education **for** social justice and sustainability,
3. and mathematics education **as** social justice and sustainability.

Our description of these approaches largely follows the works of Renert (2011), Sterling (2001), and Müller (2024, 2025b).

---

[17] In this diagram, from the perspective of Pragmatic Ethics in Mathematics, green = easy to process; yellow = can be processed with some effort; orange = more difficult to process; outside of all circles = very difficult to process. Note that the other unlabeled (blue) discursive clouds will be discussed further in the context of Figure 15.

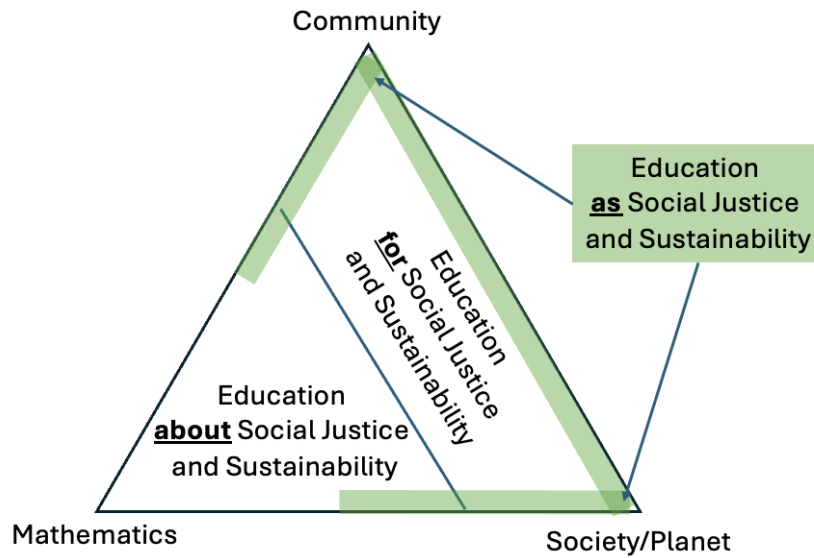

Figure 8: Mathematics education *about*, *for*, and *as* social justice and sustainability

**Mathematics Education *about* Social Justice and Sustainability**

Within mathematics education *about* social justice and sustainability, educators use mathematics to teach *about* sustainability and social justice concerns. While directed at external, local, or global ethics or sustainability concerns, a strong focus is still placed on teaching the "right kind" of mathematics. These approaches can typically be integrated into conventional mathematics learning environments without requiring significant changes in classroom culture or pedagogy. Approaches that use mathematics to teach *about* social justice and sustainability often have an intimate connection to the "Mathematics" vertex (cf. Renert, 2011; Müller, 2024, 2025b; Sterling, 2001). At the same time, they show a strong spread in the direction of both "Community" and "Society/Planet". This places such an approach to ethical and sustainable mathematics education covering the lower-left area of the triangle (Figure 8).

As an example, consider Friedman's (2022, p. 304) exercise on exponential models to model the death rates from COVID-19 for African American and White Californians, and compare them:

White Californians: $y = 0.0001444 \times (1.08991)^x$

African American Californians: $y = 0.000692782 \times (1.07156)^x$

Here, Friedman used simple mathematics to introduce the concept of environmental justice to students, writing about the experience that, "when asked to comment on what they thought might be some of the causes of the racial differences in COVID death rates, students offered suggestions such as differential access to health care, differing income levels, and other reasonable contributing factors. The instructor used this opportunity to introduce the concept of environmental justice. For all or most students, the concept of environmental justice as a type of social justice was new." (ibid., p. 304).

**Mathematics Education *for* Social Justice and Sustainability**

A more transformative approach is given by mathematics education *for* social justice and sustainable development. Here, the perspective shifts to view mathematics as a tool for positive change. It is no longer enough to teach about issues; individual and societal changes are the goal. Social justice and sustainability are no longer merely treated as a topic of (applied) mathematics. A strong focus is put on empowering students with the necessary critical and mathematical thinking skills to be advocates and agents for change. This approach, therefore, often necessitates larger pedagogical shifts and places a much stronger emphasis on both the community (e.g., the agency of students) and societal and planetary concerns.

Much of the literature advocating *for* sustainability and social justice builds on critical mathematics education, emphasising the role of mathematics in the three-step approach of understanding, critiquing, and challenging existing injustices, inequalities, and power structures (cf. Müller, 2024, 2025a, 2025b; Barwell & Hauge, 2021). In short, as Gutiérrez (2009) describes, the goal is not just to learn the game of mathematics (as depicted by classical approaches) but to change the game if necessary. Students should be empowered to have the agency to affect change in their lives and society at large. Mathematics classrooms are no longer understood as neutral learning spaces but rather as rooms governed by values, power dynamics, and other factors. The intent is to foster critical thinking, socially and morally reflective social action, and a broader critique of the role of mathematics in today's societies (e.g., Gromlich, 2021, 2024; Skovsmose, 1994a, 1994b, 2020, 2021a) and crises (Skovsmose, 2021b). This approach often questions the foundations of classical mathematics education and classical philosophies of mathematics (Müller, 2024, 2025a) and can be openly value-laden by advocating for specific changes. Thus, it spans most of the ethical concerns between "community" and "society" (as depicted in Figure 8).

As an example, consider the problem proposed by Bukoski and Erbes (2024, p. 12), and slightly adjusted here:

> An industry is being charged by the Environmental Protection Agency (EPA) with dumping unacceptable levels of toxic pollutants in a lake. Over a several month period, an engineering firm takes daily measurements of the rate at which pollutants are being discharged into the lake.

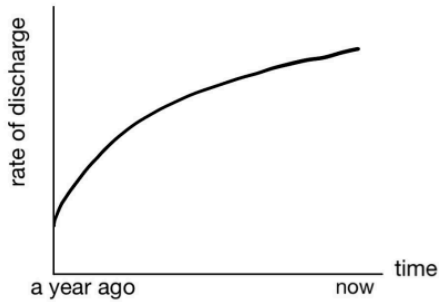 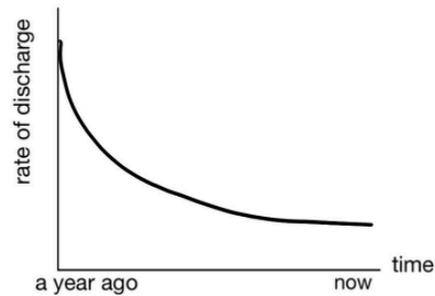

1. Suppose the engineers produce a graph similar to the graph on the left.
   a. What argument would the EPA make in court against the industry?
   b. How could the industry defend against that argument?
2. Suppose instead that the engineers produce a graph similar to the graph on the right.
   a. What argument would the EPA make in court against the industry?
   b. How could the industry defend against that argument?
3. Suppose you were hired by the EPA or the industry as a mathematical consultant. For each of the four arguments that you listed above, decide whether it would be ethical for you as a mathematician to present the argument in court. Justify your answers briefly.

This problem moves beyond merely teaching students about sustainability concerns. It asks students to make a decision about what mathematical arguments are ethical to present, and thus to use mathematics to foster a more sustainable world. However, it does not ask the students to fully reconfigure the relationship between mathematics, ethics, and sustainability, as the stronger approach of Mathematics Education *as* Social Justice and Sustainability would do.

For further mathematical exercises regarding education *about* and *for* sustainability, which are aimed at students at the qualification level II, we refer the reader to Meyer et al. (2025), who studied the effect of different exercises on the self-efficacy of students.

## Mathematics Education *as* Social Justice and Sustainability

The most holistic perspective is given by mathematics education *as* social justice and sustainability. In an ideal situation, this perspective requires all three vertices of the triangle to be re-envisioned. Mathematics is no longer seen as a neutral and purely abstract discipline, but rather as a value-laden, human, and culturally situated enterprise. The community emphasises the creation of classrooms that are collaborative, equity-focused, and built on principles of care and participatory concerns. Finally, the societal and global dimension is defined by a commitment to ethical and ecological awareness and agency-driven transformative action. Consider, for example, Valero's (2025, p. 138) analysis that "[t]o truly go in the direction of addressing the ecological crises by means [of] — or perhaps despite — mathematics education, a radical rethinking thus seems to be required: one that moves beyond the usual, learned responses to reform mathematics curricula and practices, which may easily turn into symbolic gestures not challenging the core principles and structures of the current educational paradigm."

As an example, consider the following problem from an "Introduction to Proofs" (Chiodo et al., 2025, p. 13):

> (a) Consider the statement: 'If you don't do it, then someone else will.' Express this statement in symbolic notation. Find its contrapositive and its negation, giving each both in symbols and in words.
> (b) Your boss has given you a task. The task is well within your technical capability, but you are not sure whether it would be legal or ethical. You feel uneasy, but your boss tells you: 'If you don't do it, then someone else will.' Do you think that the boss' argument is cogent? Does that depend on who you are or what the project is? How would you answer your boss?

This problem simultaneously teaches students about propositional logic (i.e., the "Mathematics" concern), about the responsibilities of being a mathematically-educated expert within a community of people (i.e., a "Community" concern), and about wider aspects of management and work in a Capitalist society (i.e., a "Society" concern). The question moves beyond asking the student to apply mathematics to understand and improve a specific situation (i.e., mathematics about and for social justice) by asking them to reconsider how mathematics should be practised and how mathematicians should respond when under ethical pressure. While the problem proposed by Bukoski and Erbes (2024) still saw mathematics as a necessary part of the solution (it was just a matter of picking the ethically good mathematics/graph), this question goes further by (implicitly) asking students about specific dynamics within mathematics ("solve this problem") and the workplace ("If you don't do it, someone else will."). Thereby, the question asks students to more fundamentally reconsider how mathematics is practiced and perceived within society.

To the authors, it seems that mathematics education *as* sustainability or as social justice can often originate from deep concerns about the community, planet, or society, which has led us to depict it as emerging from these two vertices and spreading towards mathematics. In practice, that means that literature can (implicitly) self-identify or be identified as "education *as* …" even when it puts a stronger emphasis on specific concerns. Thus, in Figure 8, we depicted such an educational approach as a green area emerging from the "Community" and "Society/Planet" vertices.

The distinction between education *about*, *for*, and *as* sustainability has strong parallels with Rieckmann's perspective on sustainable education (e.g., Rieckmann, 2016, 2017, 2018, 2021a, 2021b; Rieckmann & Holz, 2017). Rieckmann (2016) differentiates between "instrumental" ("ESD1") and "emancipatory" ("ESD2") approaches to sustainability education. ESD1 aligns strongly with mathematics education *about* sustainability, and prescription-oriented versions of education *for* sustainability. Here, the goal is to encourage specific (sustainable) actions in the students. In contrast, ESD2, characterised by deeper critical thinking, the ability to navigate complexity, and find and formulate one's own solutions to sustainability problems, is a core aspect of education *for* and *as* sustainability. Rieckmann's overall perspective is strongly competency-based to empower learners, requiring action-oriented and transformative pedagogical approaches to education. However, mathematics (education) *as* social justice and sustainability (as understood in this paper and in Müller (2025b)) goes beyond this by rethinking what *ethical and sustainable mathematics* is, and should be, beyond its education.

In his work on general education within mathematics, Winter (1996, p. 35) outlines three foundational experiences that students should have. These are:
1) "To perceive and understand phenomena of nature, society, and culture that concern or should concern us all, in a specific way."
2) "To come to know and comprehend mathematical objects and facts, represented in language, symbols, images, and formulas, as intellectual creations and as a deductively ordered world of their own kind."
3) "To acquire problem-solving skills (heuristic abilities) that transcend mathematics through engagement with mathematical tasks."

The first two experiences are essential across all educational levels of engagement: education *about*, *for*, and *as*. The third experience, however, likely assumes a particularly crucial role within the context of education *for* and *as* social justice and sustainability. In this pedagogical approach, students are not merely taught about such problems. Instead, they are explicitly encouraged to solve these problems with mathematics or potentially even without mathematics if the standard way of doing mathematics could be harmful (i.e., the limits of mathematics can be actively explored (e.g., Abtahi, 2022; Chiodo et al, 2025)). This process inherently involves an increasing critique of the role and structure of mathematics itself. Winter (1996, p. 36) elaborates that the first experience, in particular, frames mathematics as a "useful, applicable discipline" of "almost universal scope." It is precisely this aspect, i.e., the utility and perceived universality of mathematics, that also becomes an object of critical scrutiny within an educational philosophy focused on education *for* and *as* social justice and sustainability.

Winter (1996, p. 38; translated from German) writes that, "trivially, all human actions are fundamentally fallible. What is special about mathematics is that errors and misunderstandings can be objectively identified and criticised and do not represent something that a layman must believe from an expert. What is even more important for general education is that errors, misunderstandings and breaks can, by revealing their genesis, become the starting point for deeper understanding and can, so to speak, be turned into something productive." However, the idea that mathematics itself may not be the best tool to solve a problem in a socially-just or sustainable fashion (e.g., Chiodo & Müller, 2025), as is required for (advanced) education *as* sustainability, is not present in Winter's older description. This suggests that the experiences of students within socially-just and sustainable mathematics education ought to go beyond what Winter describes for general education.

Currently, most studies focus on education *about* and *for* sustainability and social justice (cf. Budikusama et al., 2024; Li, 2025; Vásquez et al., 2023; Renert, 2011; Müller, 2025a, 2025b; Makramalla et al., 2025). As explained in Müller (2025b), this has to do with the different levels of ethical engagement and awareness that mathematicians and educators can present; teaching *about*, *for*, and *as* represent first, second and third-order learning respectively (Sterling, 2001), whereby the latter usually require higher levels of ethical awareness/engagement, moving from merely seeing the issues, towards speaking about them, taking (political) action, and ultimately speaking out against injustices.

From a systems theoretic perspective, the difference between "*about*", "*for*", and "*as*" is linked to how strongly environmental pressures act on mathematics and its education. Treating sustainability merely as another application area of mathematics (education about sustainability) leaves disciplinary structures and sociomathematical norms largely intact. The environment is observed, mathematicised, and then used as content for instruction or for other mathematical work, but the autopoietic aspects of the social system of education and the system of mathematics remain largely unchanged. However, advocating for mathematics education *for* social justice and sustainability directly questions the presumed neutrality of the field and urges the discipline to pursue explicit normative goals, representing a stronger irritation that forces both the scientific and educational systems to adapt. Finally, education *as* sustainability and social justice is fully transformative, representing such an intense and persistent irritation that it forces a re-organisation of the structural couplings between all systems, and potentially even a change in the fundamental differentiation between them. The boundaries between the scientific system, the educational system, and their environment become blurred as they are forced to co-evolve into a new, more integrated structure.

## Archetypal Educators

To further nuance the location effect, and to understand some concrete educational positions behind "mathematics education about, for, and as …", we will now consider the five educator archetypes which Ernest (1991) identified in 20th-century British mathematics education. These are represented in Figure 9. Müller (2025a) argues that these archetypes can also be found in discourses on ethics in mathematics and its education. Over time, they have developed educational goals that put them into (direct) conflict, with each group fighting for their vision of education. Ernest (2019a, p. 86; 2019b, p. 72) argues that a new balance may be necessary because historically it was not that with the most ethical vision, but rather that with the most power, dominance and utility which has won. Understanding what archetype someone embodies and where they are located will help to provide a more nuanced understanding of specific sub-discourses later on. From a systems-theoretic perspective, these archetypes represent distinct patterns of structural coupling between the social communication system of education with its environment and the system of mathematics.

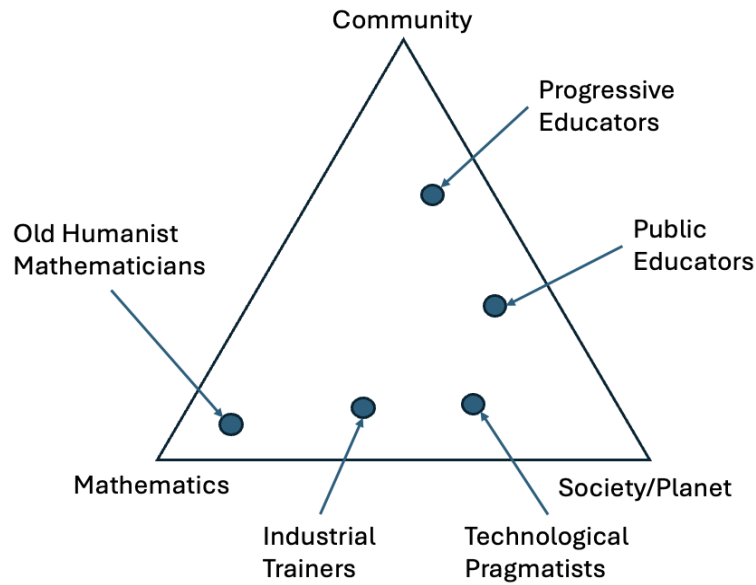

Figure 9: Archetypal educators

## Industrial Trainers and Technological Pragmatists

The archetypal industrial trainer sees mathematics as a utilitarian tool to equip learners with the basic skills required for employment (Ernest, 1991, p. 147). As their teaching prioritises perspectives on job readiness, their ethical considerations may centre on issues of professional conduct, diligence, and the responsible application of mathematical skills. They are most closely aligned with educational paradigms of teaching about ethics and sustainability because these require the least socio-cultural and economic adjustments. Their emphasis would be strong on the "Mathematics" and "Society/Planet" vertices – albeit potentially quite unpolitical. Ernest (2019a, p. 84; 2019b, p. 70) argues that "[t]his group does not want education politicised in order to prepare a demanding and non compliant workforce", but rather wants a "social training in obedience."

Technological pragmatists also view mathematics from a utilitarian perspective, understanding it as a powerful tool for practical problem-solving and technological innovation. As Ernest (1991, pp. 162 - 163) notes, their vision of (school) mathematics contains two parts: "First, there are the pure mathematical skills, procedures, facts and knowledge. These are the dry bones of the subject, which are simply tools to be mastered. Second, there are applications and uses of mathematics. This is the vital, living part of mathematics, which justifies and motivates the study of the subject." Their teaching strongly values mathematics to efficiently solve problems, and their ethical discussions may centre on the responsible use of technology and the broader societal implications of technological advancements. This orientation also places them near the "Mathematics" and "Society/Planet" edge of the triangle, albeit slightly further in the direction of "Society/Planet", as they can be more pragmatic about specific mathematical concerns so long as it justifies or relates to higher external goals. Ernest (2019a, p. 85; 2019b, p. 71) sees in this group a "meritocratic vision of society" enabling individuals to climb the socio-economic ladder through education.

From a systems theoretic perspective, both industrial trainers and technological pragmatists represent a strong coupling of mathematics with its economic environment.

### Old Humanist Mathematicians

This archetype strongly appreciates mathematics for its own sake. As Ernest (1991, p. 172) writes, "The central element of this ideology is that education and knowledge are a good, an end in themselves, and not a means to a baser, utilitarian end." They emphasise the rigour and elegance of mathematical proof and the internal consistency of mathematics, as well as mathematics' timeless and universal nature (Ernest, 2019a, 2019b). This position aligns strongly with the "Mathematics" vertex, with some consideration for the "Community" of mathematicians, students, and educators who preserve and advance this (cultural) heritage, as well as some "societal" aspects because mathematics is generally seen as a tool of progress and this therefore furthers the case for more mathematical study (cf. Ernest, 2019a, p. 85; Ernest 2019b, p. 71). From a systems theoretic perspective, the old humanist mathematicians represent a strong coupling to the scientific system. Old humanist mathematicians particularly care about the autopoiesis of mathematical knowledge, and have for a long time experienced few external irritations. For example, Hersh wrote about the development of the American Mathematical Society's Code of Ethics in 1995[18]:

> "Almost 20 years ago, I noticed that the American Mathematical Society, unlike most professional organizations such as chemists, geologists, statisticians, etc., had no official Code of Ethics. This observation stimulated some pondering on my part, which was vented in a talk to a regional meeting of the Mathematical Association of America, and an article in the Mathematical Intelligencer in 1990. I conjectured that part of the reason such an official code was absent was that, while research in pure mathematics can be good or bad, it is pretty harmless, compared to the possible danger of unethical behavior in, say, nuclear physics or chemical engineering." (Hersh, [R], 2008, p. 1)

---

[18] The ethical aspects of pure mathematics have become clearer since then (Chiodo & Clifton, 2019).

In this sense, the old humanist mathematicians, as described by Ernest (1991, p. 178), "view [mathematics] as pure and unrelated to social issues, so no room is allowed for the accommodation of social diversity. Mathematics is objective, and attempts to humanise it for educational purposes, however well intentioned, compromise its essential nature and purity." Within ethical and sustainable mathematics and its education, these positions may be found in those who care more about developing the right mathematics than solving the right problem (cf. Chiodo & Müller, 2025). Teaching is primarily about teaching mathematics, and social justice and sustainability are often restricted to being areas of mathematical applications. From this perspective, deeper engagement with both ethics and sustainability, which also attempts to reform how mathematics is done and taught, may be interpreted as a "joke" rather than being taken seriously (cf. Shah, [R.], 2024; 2025).

### Progressive Educators

This archetype focuses on student-centred learning and individual development, among other things, by prioritising supportive, equitable, and engaging learning situations. Progressive educators depict a certain form of "liberal romanticism" (Richards, 1984, as quoted in Ernest, 1991, p. 185). At the core is the learner, who needs to be protected and nurtured in their growth through an appropriate educational ideology (Ernest, 1991, pp. 185 - 188). They understand that merely using mathematics to teach about sustainability is likely insufficient to establish this learning environment. As Ernest (2019a, p. 85; 2019b, p. 71) argues, this group of educators wants students to be creative and develop robust forms of self-expression and confidence that make the individual flourish, albeit, as he argues, potentially at the expense of other forms of social awareness and the promotion of social goods. This places progressive educators strongly within the "Community" vertex, as their concerns are, first and foremost, centred around the development of learners. This can be understood as education *for* social justice and sustainability, but also as localised education *as* social justice and sustainability.

From a systems-theoretic perspective, progressive educators care for the autopoiesis of the educational system and the psychic systems of the learners. They focus on communicative patterns (pedagogies) that are internally oriented towards fostering specific cognitive operations within learners.

### Public Educators

Viewing mathematics as crucial for empowering citizens to participate in democracy, this archetype focuses on developing mathematical literacy for socio-political awareness and action (Ernest, 2019a, pp. 85 - 85; Ernest, 2019b, p. 71). As Ernest (1991, p. 199) argues,

"[t]he goal of this position is the fulfilment of the individual's potential within the context of society. Their teaching uses mathematics as a tool to understand, critique, and change societal structures, often with explicit ethical stances concerned with access, equity, and social justice. Thus, the aim is the empowerment and liberation of the individual through education to play an active role in making his or her own destiny, and to initiate and participate in social growth and change". From this perspective, "[i]t is of democratic importance, to the individual as well as to society at large, that any citizen is provided with instruments for understanding th[e] role of mathematics [in society]. Anyone not in possession of such instruments becomes a 'victim' of societal processes in which mathematics is a component" (Niss, 1983, p. 248, as quoted in Ernest, 1991, p. 205). Ernest (2019a, p. 85; 2019b, p. 71) sees in this archetype a group of educators explicitly trying to politicise education, and that the individual well-being of students may be traded in for larger social goals.

This locates public educators with a strong focus on both the "Community" and "Society/Planet". Following Ernest's description, they are placed further in the societal corner than progressive educators. Both archetypal educators are likely involved in teaching *for* and *as* sustainability. Together with progressive educators, public educators may display high levels of ethical and sustainable awareness as defined by Rycroft-Smith et al. (2024) and Chiodo & Bursill-Hall (2018). From a systems theoretic perspective, they represent a strong coupling with the socio-political environment.

### Academic and Student Trainers

However, the mapping of Ernest's (1991) archetypes reveals a notable gap in the triangle between the "Mathematics" and "Community" vertices. The Old Humanist Mathematicians' educational goal is strongly focused on the transmission of cultural heritage, and as described by Ernest (1991, 2019a, 2019b), they typically do not have strong community concerns. This suggests the possibility of additional educator archetypes whose concerns are strongly rooted in both areas, without significant orientation towards the "Society/Planet" vertex. We posit two such archetypes: the Academic Trainer and the Student Trainer. Both share a common viewpoint of "mathematics for the sake of mathematics," but their primary motivations diverge in a subtle but important manner.

The Academic Trainer views education as a means to an end: the preservation and advancement of research mathematics. Potentially found in the setting of research universities, their focus is on producing the next generation of research mathematicians. For them, educational practices are designed and evaluated based on their effectiveness in identifying and cultivating future researchers. Their mantra is "quality, not quantity," where quality is measured by the potential to produce research-level mathematics. A successful student, from this perspective, is one who ultimately pursues a career in academic mathematics. While they care deeply about the community of students and educators, this concern is instrumental to the larger goal of sustaining the body of mathematical knowledge. They perceive mathematics as a (static) body of knowledge to be maintained and expanded where possible. Consider, for example, Körner's (n.d.) "Unofficial Guide to Part III" of the Mathematical Tripos at Cambridge, who argues to study in a way that (from his perspective) most likely produces a high-quality research mathematician, and who argues that the course itself is only designed for those of sufficiently high mathematical calibre: "you should be able to profit from Part III if you are in the top 10% of mathematicians graduating in your country and you are prepared to work very hard" (ibid., p. 3).

The Student Trainer, in contrast, views mathematics as a vehicle to further education. Their primary objective is to have as many students as possible engage with and appreciate mathematics as a fulfilling activity in its own right, akin to a sport or leisure pursuit. They value broad participation and individual growth across all levels of mathematical ability, often engaging in activities like mathematics olympiads. For student trainers, mathematics is a kinetic activity that comes alive through engagement, rather than a static body of knowledge to be passively learned. A successful student is one who develops their mathematical potential to the fullest, regardless of their ultimate career path. While they are deeply invested in mathematical research as a source of rich content for their students, their fundamental goal is the enrichment of the student community. One such example can be found in the mathematician Imre Leader, who argues that, "I think the main thing for high school students is that if they like the stuff, it almost doesn't matter what they do; it's the number of hours that matter. [...] The number of hours you spend, say, in math, is just good for your math brain. It makes you more mature, and it almost doesn't matter what you study. So just spend time (doing what you like)" (Zhong & Nguyen, 2024).

Though located closely together on the ESCT (Figure 10), these two archetypes exhibit a key distinction in their orientation. The Academic Trainer looks toward the "Mathematics" vertex, using the community to serve its advancement, whereas the Student Trainer looks toward the "Community" vertex, using mathematical knowledge to serve its development. They are ideologically aligned in their shared focus but are, in effect, looking in opposite directions. This comes from the difference in how they view mathematics, with one seeing it as a static object whose origins are somewhat inconsequential, and the other seeing it as a kinetic object which only exists when it is being actively carried out. While the Old Humanist

Mathematicians can have a superficially similar vision of teaching the future mathematical elite (Ernest, 2019a, p. 85; 2019b, p. 71), the academic and student trainers can be more utilitarian, more social, and less pure in their educational approach – in short, there is now a place for the humanisation of mathematics and for political educational concerns *if* they help to better teach advanced mathematics.

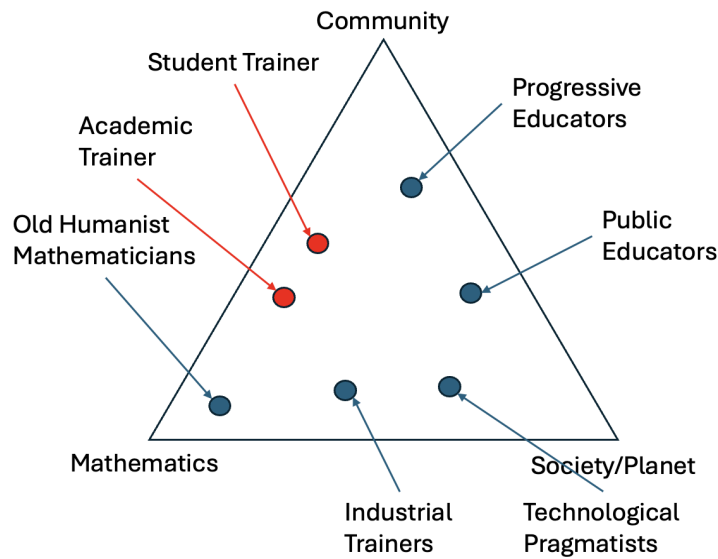

Figure 10: Mathematical trainers

## **The Magical Educator**

Having established two new educational archetypes, the ESCT still reveals a curious void in its centre when mapping various educator archetypes. This empty space suggests the absence of a "magical educator" who perfectly balances mathematical, community/educational, and socio-planetary concerns (see also Abtahi's (2022) discussion of "What if I was harmful?" in the context of teaching mathematics). We argue that this is not an accidental omission but rather a significant finding of the framework. It indicates that this ideal, balanced position is not a practical reality within current academic and educational contexts. The archetypes tend to cluster closer to the vertices or along the edges, suggesting that educators are typically driven by a primary concern or a combination of two, rather than a harmonious integration of all three.

A systems-theoretic analysis explains why this central void may exist. The "Mathematics" and "Community" vertices can be reframed as autopoietic mathematical and educational systems. Each has its own internal logic, norms, and self-reinforcing practices. It is possible for an individual to find a stable and functionally successful equilibrium between these two systems. For example, an educator focused on training future mathematicians serves the reproductive needs of both the mathematical and educational communities. In this stable position, the need for further ethical engagement can diminish because the existing arrangement is perceived as successful on its own terms. The "Society/Planet" vertex, in contrast, is not a self-contained system but rather the shared environment from which persistent "irritations," such as social injustice or climate change, emerge, and which strongly values specific types of (economic) education, as evidenced by the Industrial Trainers and the Technological Pragmatists.

This dynamic has profound consequences for engagement with broader ethical issues. When a stable equilibrium is achieved between the mathematical and educational systems, the troubling irritations from the socio-planetary environment can be effectively filtered out as noise. Consequently, meaningful engagement with socio-planetary concerns often requires approaching it from a strong grounding in one of the other systems, and thus hugging an edge of the triangle. Educators positioned in a stable balance between mathematical and community concerns may not feel the systemic pressure or need to engage with the complex and often disruptive issues of ethics and sustainability that define the third vertex.

This difficulty is not unique to the "Mathematics-Community" edge. A similar dynamic can be observed along the other edges of the triangle, although the reason may be subtly different. An educator situated between "Society/Planet" and "Mathematics," for instance, may not be in a state of stable equilibrium but rather in a constant state of tension, perpetually engaged in mediating between powerful societal demands and the traditions of the mathematical system (cf. Figure 6). This continuous intellectual pull between two strong vertices leaves comparatively little room to fully engage with the third, making a move toward the balanced centre of the triangle equally difficult. Similarly, someone situated between "Community" and "Society/Planet" may be constantly experiencing tensions between individual students' needs and socio-planetary demands, making it difficult for him to relocate to a central position inside the triangle (cf. Figure 6). In this context, Meyer (in preparation) observed multiple general tensions, including 1) factual vs normative, 2) unambiguity (mathematics) vs ambiguity (in sustainability), 3) what topic to teach next, 4) if it should be from mathematics or sustainability, 5) balancing mathematical, STEM, and ethical competencies, and 6) how to sell negative aspects positively enough to keep students engaged. These tensions translate into tasks, or challenges, for each archetypal educator, and different educators will (likely) answer them differently.

## School and Higher-Mathematics

These archetypal educators underline that the field has a porous boundary between ethical and sustainable concerns about school and university (as similarly discussed in Müller, 2025a).[19] However, despite a porous boundary of concerns, we can observe that, overall, the literature which focuses on higher-mathematics is much more aligned with concerns about society, the planet, and mathematics, while the literature which focuses on school mathematics has a much stronger connection to the "Community" vertex and societal issues (as depicted in Figure 11, and further explored in the context of Figure 15). This goes so far that Borovik (2025, p. 1) "feel[s] that real mathematics and real mathematicians are disappearing from the picture" because so much of the (mathematics (education)) literature on ethics does not strongly focus on the development of mathematics.

From a systems-theoretic perspective, the school versus university shading in Figure 11 illustrates how different subsystems exhibit different patterns of structural coupling. The system of university mathematics is strongly coupled to the scientific system's (mathematical) research frontiers, processing irritations primarily from within the discipline. In contrast, the social system of school mathematics is more strongly coupled to the psychic systems of students and the broader socio-political environment, thus processing a wider range of social and personal irritations (cf. Watson, 2023).

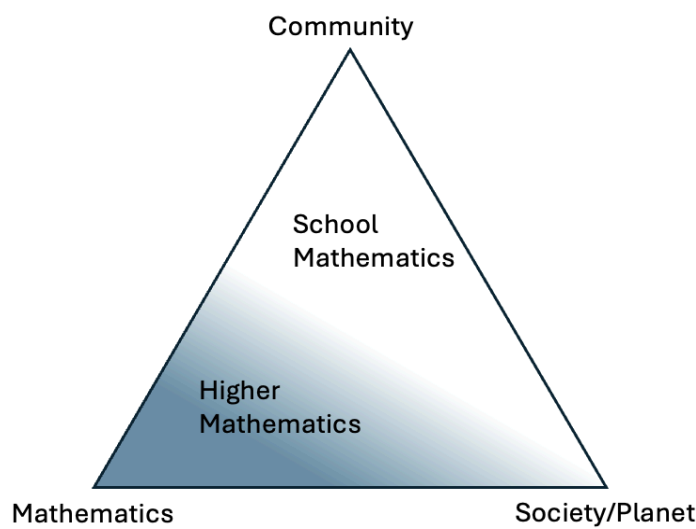

Figure 11: Higher versus school mathematics

---

[19] For example, the levels of ethical engagement developed in the context of university mathematics (Chiodo & Bursill-Hall, 2018) were then transferred by Rycroft-Smith et al. (2024) into the school context.

When focusing on the school level, the discourses on ethics and sustainability in mathematics and its education frequently centre on the immediate learning environment and the individuals within it. This includes the nature of the teacher-student relationship, student well-being, gentle and mindful teaching, the establishment of inclusive and equitable classroom cultures, and how the presentation of mathematical content can impact understanding, identity, and self-expression (cf. Müller, 2025a; Müller, 2025b; Rycroft-Smith et al., 2024). Ethical and sustainable concerns are typically integrated through age-appropriate topics (cf. Budikusama et al., 2024; Müller 2025b; Li, 2025; Wilhelm, 2024), and care must be taken to consider how their framing affects students (e.g., Meyer et al., 2025). The main focus generally lies on teacher education (e.g., Rycroft-Smith et al., 2025, Vásquez et al., 2023), the students and their development, with only a secondary focus on the development of (research) mathematics and wider society.

At the university level, the focus tends to shift towards the ethics of mathematical practice itself, particularly concerning research and advanced mathematical applications (Chiodo & Müller, 2024; Müller et al., 2022). While community concerns, such as equity, student well-being, diversity, and inclusion, are becoming increasingly important in higher mathematics education (e.g., Buell & Piercey, 2022), many discussions still revolve around the professional responsibilities of mathematicians (e.g., Tractenberg et al., 2024), the debated neutrality of pure mathematics (e.g., Ernest, 2024b), the potential for the misuse of mathematical work such as in data science, artificial intelligence, climate change and financial modelling, and gerrymandering (e.g., Thompson, 2022), and the broader socio-ecological impacts of mathematical research and innovation (e.g., Müller, 2025a; Chiodo & Vyas, 2019; Chiodo & Clifton, 2019; Müller et al., 2022). It typically involves applications, e.g., the modelling of global environmental challenges and research into sustainable technologies (e.g., Hersh [, M.], 2006), or the analyses of mathematics's role in shaping socio-economic and political paradigms (e.g., Müller & Chiodo, 2023). Next to the "Mathematics" vertex, the "Society/Planet" vertex is particularly strong here, with applied mathematics and modelling being closely linked to it.

Overall, it seems to the authors that the origins of ethical concerns regarding "school education" are much more firmly rooted in community concerns, while their origins in "higher education" are much more rooted in mathematical concerns (Figure 12), with notable exceptions when educators deviate from this expectation and, essentially, focus on other vertices, as also explored in the Academic Trainer and Student Trainer archetypes of the previous section. The white gap and the arrows depicted in Figure 12 represent both how different these concerns can be, but also that it is possible for (individual) educators to cross this divide and to move beyond social (and sociomathematical) norms and expectations.

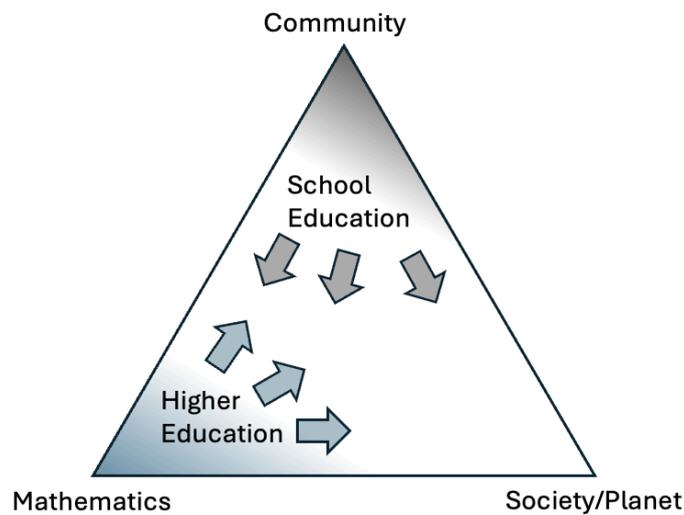

Figure 12: Origins of ethical concerns in higher and school education

As a concrete example from higher education, consider the British statistician Jane Hutton, whose career exemplifies a clear trajectory from foundational quantitative statistical science towards its application in complex societal domains like medicine, ethics, and law. With a PhD in statistics (Mathematics Genealogy Project, n.d.), her early academic appointments and publications focused on medical statistics and statistical methods (e.g., Ashby et al., 1993; Carr et al., 1990; Escala et al., 1989; Hartley et al., 1988; Hutton, 1990; Hutton & Owens, 1993; Modi & Hutton, 1990a, 1990b; Moore et al., 1991; Parys et al., 1989, 1990a, 1990b; Smith et al., 1991). This engagement in high-stakes medical contexts and her upbringing under apartheid in South Africa appear to have spurred a deeper inquiry into the principles governing such research (Hutton, 2018). This is also evidenced by a continued stream of publications on ethics and professional issues since the mid-1990s, when she began writing on how statistics is essential for professional ethics and by a series of philosophical articles exploring the ethics of randomised controlled trials, research in developing countries, and other ethical aspects of the statistical work itself (e.g., Ashcroft et al., 1997; Buyse et al., 1999; Hutton, 1995, 1996, 1998, 2000; Hutton & Ashcroft, 1998).

Her ethical grounding provided a natural bridge to the legal field, where she has served as an expert witness in over 500 cases since 1995 (Hutton, 2025, p. 2), before eventually becoming a founding member of the Royal Statistical Society's "Statistics and the Law" section (University of Warwick, n.d.), and bringing ethical and legal insights back into mathematics, statistics, and law education itself (Hutton, 2025, p. 5). She said about her experience that "quite a few of my colleagues don't think I'm a real mathematician, but quite a few of my other colleagues think that I'm far too theoretical. I'm always [...] somebody who's on a spectrum, and actually I find my niche trying to translate from one community to another or between various communities" (Hutton, 2023).

For an example from school education, consider the American philosopher, civil rights activist, and mathematics educator Robert (Bob) Moses. After studying philosophy and French, he began teaching mathematics (Levenson & Medina, 2021). Inspired and challenged by the picketing and sit-ins at lunch counters of Black people in the South of the United States, he became a civil rights organiser in the state of Mississippi, recalling his initial motivation in his book *Radical Equations: Civil Rights from Mississippi to the Algebra Project*:

> "The sit-ines woke me up. Until then, my Black life was conflicted. I was a twenty-six-year-old teacher at Horace Mann, an elite private school in the Bronx, moving back and forth between the sharply contrasting worlds of Hamilton College, Harvard University, Horace Mann, and Harlem. The sit-ins hit me powerfully, in the soul as well as the brain. I was mesmerized by the pictures I saw almost every day on the front pages of the New York Times – young committed Black faces seated at lunch counters or picketing, directly and with great dignity, challenging white supremacy in the South. They looked like I felt."
> (Moses & Cobb, 2002, p. 3)

Out of the injustices he experienced within his own Black community, the classroom, and society at large, eventually grew *The Algebra Project,* which put as its mission "that just as people at the bottom levels of society in Mississippi require access to the right to vote to gain political power and access to full U.S. citizenship, young people need access to algebra to gain full citizenship in the 21st Century" (The Algebra Project, 2025).

## Neutrality, Purity, Universality (NPU) Assumptions

Having mapped out differences between ethics and sustainability in the discourses surrounding school and higher mathematics and different archetypal educators, we can now consider deeper philosophical assumptions about the nature of mathematics itself. As explored in Müller (2025a), the Neutrality, Purity, Universality (NPU) assumptions are a cluster of values and beliefs[20] that shape the perception of the nature of mathematics:

- Neutrality describes the conviction that mathematics is apolitical, value-free, objective, and devoid of human interests and societal influence.
- Purity has a dual meaning of "pure mathematics" (as a field separate from applied mathematics) and as a normative position (moral purism). It can view external input as a potential form of contamination, and sees mathematics as an enterprise that (almost always) does good.
- Universality asserts that mathematical truths transcend cultures and different times. In this perspective, mathematics is more likely to be discovered than invented.

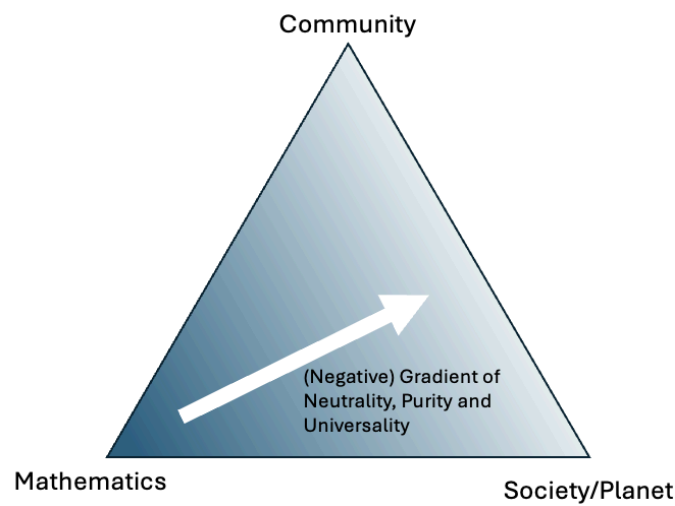

Figure 13: The NPU gradient – NPU decreases in the direction of the white arrow.

These values and beliefs are typically associated with classical perspectives on mathematics (e.g., Platonism, formalism, logicism), and alternative perspectives (e.g., social constructivism) move away from them, as do those scholars who deeply value community and societal concerns (cf. Müller, 2025a). In Figure 13, we visualise this through the shading and the arrow pointing away from the "Mathematics" vertex, depicting how the NPU decreases with distance to the "Mathematics" vertex. It is worth noting that the shading between school mathematics and higher mathematics may actually be different because ethical and sustainable school mathematics already tends to be skewed towards the "Community" vertex, while ethical and sustainable (university) mathematics is often more skewed towards research and the applications of mathematics. Hence, within higher mathematics, it is more likely to find NPU assumptions, even among those who are concerned about societal issues (cf. Ernest, 2020a, 2020b, 2024b; Müller, 2025a; Müller et al., 2022).

---

[20] For further reading and an overview on values and beliefs related to ethics and mathematics, we refer the reader to Ernest (e.g., 2016, 2020a), noting here in particular that "believing there is no ethics in mathematics" is a common theme found among mathematicians (Chiodo & Bursill-Hall, 2018, p. 5) and mathematics educators (Rycroft-Smith et al., 2024).

From a systems perspective, autopoiesis stands at the centre here. Claims of neutrality, purity, and universality are mechanisms through which mathematics reproduces its structures and asserts independence from its environment. Narratives that dispute these claims by insisting that mathematics is value-laden, culturally situated, and historically contingent, appear as irritations that threaten disciplinary self-maintenance – or are simply ignored as noise, meaning that mathematicians may favour those (philosophical) beliefs which support the mathematical system's autopoiesis, potentially foster exclusion, elitism, hierarchical forms of knowledge, and ethical indifference (Spindler, 2022). In other words, different philosophical beliefs can be "untenable" and rather "unsympathetic" for mathematicians, as Hardy explained:

> "It is often said that mathematics can be fitted on to any philosophy, and up to a point it is obviously true. Relativity does not (whatever Eddington may say) compel us to be idealists.' The theory of numbers does not commit us to any particular view of the nature of truth. However that may be, there is no doubt that mathematics does create very strong philosophical prejudices, and that the tests which a philosophy must satisfy before a mathematician will look at it are likely to be very different from those imposed by a biologist or a theologian. I am sure that my own philosophical prejudices are as strong as my philosophical knowledge is scanty. One may divide philosophies into sympathetic and unsympathetic, those in which we should like to believe and those which we instinctively hate, and into tenable and untenable, those in which it is possible to believe and those in which it is not. To me, for example, and I imagine to most mathematicians, Behaviourism and Pragmatism are both unsympathetic and untenable. The philosophy of Mr. Bradley may be just tenable, but it is highly unsympathetic. The Cambridge New Realism, in its cruder forms, is very sympathetic, but I am afraid that, in the forms in which I like it best, it may be hardly tenable. 'Thin' philosophies, if I may adopt the expressive classification of William James, are generally sympathetic to me, and 'thick' ones unsympathetic. The problem is to find a philosophy which is both sympathetic and tenable; it is not reasonable to hope for any higher degree of assurance." (Hardy, 1929, p. 3)

These foundational perspectives – from the systems-theoretic interpretation and the location effect to the mapping of educator archetypes and NPU assumptions – now help us to analyse different discourses in Part II.

# Part II: Mapping of Different Ethics and Sustainability Discourses

To map the various discourses onto the ESCT, we employ a qualitative methodology inspired by the discourse analysis of Laclau and Mouffe (2020) by focusing on three key areas: articulation practices (Table 2), the uncovering of identities, hegemonies, and power relations (Table 3), and the identification of discursive nodal points (Table 4). Their approach is particularly fitting, as it was developed specifically to challenge the historical pre-determinism inherent in rigid structuralist theories like traditional Marxism. We adopt this critical lens because we recognise that the fixed, triangular structure of the ESCT itself risks imposing its own form of structural determinism, thereby replicating, for instance, the divide between "nature" and "culture" that many of the mapped discourses aim to overcome, an issue we revisit with ethnomathematics later on. Therefore, using Laclau and Mouffe's emphasis on articulation, power relations, and nodal points allows us to treat the map not as a static grid, but as a dynamic tool to reveal the very tensions its structure might otherwise conceal.

More concretely, while the ESCT provides a structuralist snapshot of concerns – especially with its "[M]athematics" vertex representing a body of formal knowledge – the lens of Laclau and Mouffe introduces a post-structuralist balance. It allows us to interpret the "location effect" with greater discursive precision, understanding that a discourse's position is continuously (re-)constructed through its articulation practices. These practices, detailed in Table 2, relate to educational contexts, schools of thought, and perspectives on the purpose of mathematics education. Different discourses adopt distinct articulative styles, often corresponding to their location on the triangle: one finds relatively political language near the "Society/Planet" vertex, comparatively rationalist articulations close to the "Mathematics" vertex, and a focus on individual existence and well-being in discourses near the "Community" vertex.

Another core concept in this analysis is antagonism. Laclau and Mouffe (2020) argue that social groups define their positions and interests not in isolation, but through the experience of conflict with opposing views. These antagonisms are essential for understanding a discourse's boundaries. In the literature on ethics and sustainability, such conflicts are rooted in fundamental power relations (e.g., different levels of ethical and political engagement), deep-seated philosophical assumptions (e.g., about the NPU of mathematics), and questions of identity (e.g., Ernest's archetypal educators). Table 3 groups these concerns.

Finally, we analyse each discourse's nodal points, i.e., the privileged signifiers that orient a discourse at a fundamental level by temporarily fixing meaning. For example, this paper establishes a nodal point by defining the "Mathematics" vertex as a formal body of knowledge. Laclau and Mouffe argue that a discourse's success depends on its ability to fix such signifiers and develop a "chain of equivalence" between them. Concrete examples include how well-being has become a prominent community concern or how data bias has become a nodal point for debates about AI. We identify these key concepts by examining which aspects of the mathematical workflow (Chiodo & Müller, 2025) a discourse prioritises, as summarised in Table 4.

| **Articulation Practices** | |
|---|---|
| **Educational Context** | Is the primary focus on the school or university level? How does this context shape the discussion? (as discussed in the context of Figure 11) |
| **Authorial Stance** | If, and how do individual papers (and thus entire discourses) make their own position on mathematics education about, for, or as social justice and sustainability clear? (as discussed in the context of Figure 8) |
| **Schools of Thought** | To which broad educational school of thought (e.g., critical mathematics education) do the papers inside a discourse belong? |
| **Expressed and Hidden Concerns** | Do the authors make their concerns clear or are they left implicit? How can their concerns be understood from the proposed triangle-based system's theoretic perspective? |

Table 2: Articulation practices

| **Uncovering Identities, Hegemonies and Power Relations** | |
|---|---|
| **Levels of Ethical Engagement:** | What level of ethical awareness do the papers (and thus the discourse) display, and where do they put it into action, if at all (adjusted from Rycroft-Smith et al., 2024)?<br><br>● Actively obstructing efforts to address ethics and sustainability in mathematics<br>● Believing there is no ethics and sustainability in mathematics<br>● Realising there are ethical and sustainable issues inherent in mathematics<br>● Doing something: speaking out to other mathematicians and educators<br>● Considering seats at the tables of power<br>● Calling out the bad mathematics of others |
| **Core Assumptions** | What assumptions about the neutrality, purity, and universality of mathematics do the papers (and thus the discourse) uphold or challenge? (as discussed in the context of Figure 13) |
| **Philosophical Foundation** | What philosophy of mathematics appears to underpin the main arguments? |
| **Educational Position (Identity)** | What type of educator do the papers (and thus the discourse) appear to assume? (as discussed in the context of Figures 9 and 10) |

Table 3: Uncovering hegemonies and power relations

| Identification of Mathematical Nodal Points | |
|---|---|
| **Workflow Engagement:** Which aspects of the responsible mathematical workflow do the papers (and thus discourses) mainly engage with (according to Chiodo & Müller, 2025 and Chiodo & Müller, forthcoming)? | Mathematical concerns: handling data and information, data manipulation and inference, mathematisation of problems, etc. |
| | Community concerns: diversity and perspectives, communication, politics of mathematics, etc. |
| | Societal and planetary concerns: large scale impact of mathematics, etc. |

Table 4: Identification of mathematical nodal points

In summary, we used the "entirety" of the literature to construct specific triangular organisational principles (the preceding triangular figures of this paper) and turned these into a qualitative discourse analytic method, which we can then apply to specific discourses themselves.[21] In doing so, we have constructed a discourse analytic lens through which we can conduct a qualitative, interpretive meta-review to locate each discourse within the triangle (Figure 14). The questions in Tables 2 - 4 strongly guided the analysis, but did not fully determine it; the authors' decade-long experience of working in the area has also influenced their judgment and the placement of individual discourses when synthesising the findings.[22]

---

[21] One way to understand these discourse clouds is as a representation of what to expect from a given discourse. Because these are qualitative, interpretive visualisations representing the discourse's apparent conceptual center of mass and its intellectual spread across the ESCT's framework's core concerns, some discourses have comparatively small/large clouds, even though the quantity of their scholarship is actually greater/smaller than some of the other approaches (cf. Müller, 2024; Müller 2025a; Dubbs, 2021).

[22] How some of the authors positioned their own pragmatic approaches to Ethics in Mathematics in relation to each question is further explored in Table 5 of Appendix A. Readers interested in applying

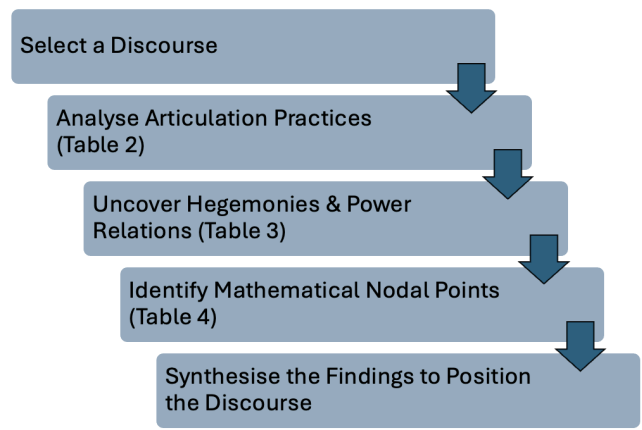

Figure 14: Workflow to map a discourse

The selection of discourses is partly due to space constraints and partly explained by the location effect of what the authors have perceived as discourses about the ethics and sustainability of mathematics and its education at large, rather than subject-specific discourses (e.g., the ethics of finance, cryptography, or AI). The discourses were mapped out in the order given in this paper. The reception of Levinas, calls for Hippocratic oaths, and the positions of professional mathematical societies were located first because some of the authors had previously written review-type articles about these. Next, pragmatic approaches to ethics and sustainability were mapped out because it is where some of the authors self-located their scholarship. Having established discourses near the vertices and (some of) the authors' positions, other discourses were mapped out. This order was understood as a rational methodology[23]: the subjectivity from some of the authors' positions becomes visible, but given that three other discourses were already mapped onto the triangle, later discourses could now be mapped in a more complex relation.

The sociopolitical (Gutiérrez, 2013) and ethical turns (Müller, 2025a; Ernest, 2024a) have significantly shifted the discourses in mathematics and its education. In Figure 15, this shift is visually represented by an arrow pointing away from the "Mathematics" vertex, indicating an overall move away from purely mathematical concerns to community and societal/planetary concerns, and their sociopolitics and ethics. Figure 15 shows the centrality of some discourses (e.g., socially just and sustainable modelling), hinting at their ability to bridge different discourses, which we detail along with all the other discourses in the remaining part

---

the framework, can also find further guidance for understanding their own position in Appendix C, Table 7; in addition to the Tables 2 - 4 above.

[23] Start with areas where the authors have written reviews; map out some of the authors' own positions; check for consistency; adjust if necessary; map out further discourses; checking for consistency each time.

of the paper. Having points of contact with many other discourses means that these discourses are rarely perceived as noise.

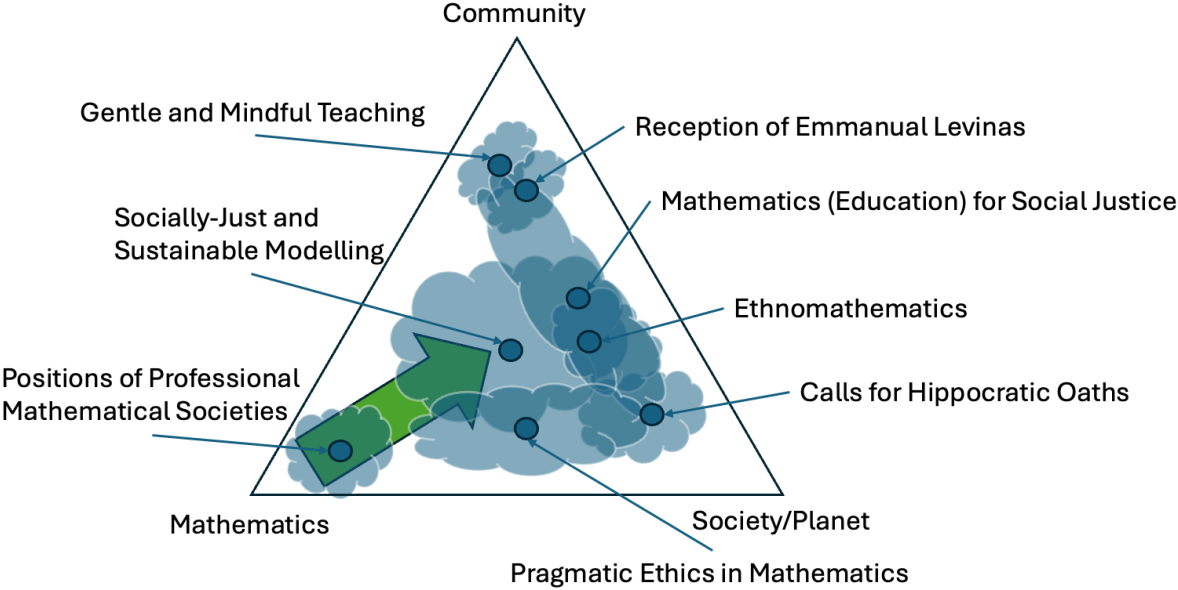

Figure 15: Discourse clouds of ethics and sustainability

As positions within this triangle illustrate how systems irritate each other, the discourses situated between "Community" and "Society/Planet" represent communications within the social system of education, processing strong irritations from the system's environment. Discourses situated between "Mathematics" and "Society/Planet", such as the pragmatic Education in Mathematics tradition developed by the Ethics in Mathematics Project (Müller, 2024, 2025a; EiMP, 2025), demonstrate the scientific system of mathematics reflecting on its impact. Discourses that are situated between the "Mathematics" and "Community" vertices represent couplings between the two systems and how these two systems can irritate each other. For example, the discourse on gentle and mindful teaching can be seen as the social system of education producing a set of communications that reacts to irritations from the rationalist, utilitarian communications of the scientific system.

Below, we briefly summarise each of the following discourses: Reception of Emmanuel Levinas, Calls for Hippocratic Oaths, Positions of Professional Mathematical Societies, Pragmatic Ethics in Mathematics, Mathematics (Education) for Social Justice, Socially-Just and Sustainable Modelling, Gentle and Mindful Teaching, and Ethnomathematics. Their order corresponds to the order of our initial analysis to assist readers wishing to replicate aspects of the mapping.

# Reception of Emmanuel Levinas

Emmanuel Levinas is an ethicist who has received much attention in discourses on ethics. As shown by Müller (2025a, pp. 31 - 36), scholars concerned with the immediate, interpersonal dynamics of teachers and learners primarily engage with his scholarship. This places the reception of Levinas by education scholars, and the consequent discourses building on it, near the "Community" vertex of the triangle. While the discourse surrounding Levinas prioritises the responsibility and relationship to the "Other" (person), there is also a non-negligible concern for the betterment of society, for example, in the works on "socially response-able" mathematics education, aimed at empowering students and teachers to engage with (i.e., respond to) the world and people's needs (Atweh, 2007, 2011, 2012, 2013a, 2013b; Atweh & Ala'i, 2012; Atweh & Brady, 2009; Atweh & Swanson, 2016). However, the nature of his reception is also rooted in a fundamental critique of mathematics as a neutral, universal, and purely rational discipline. Consider, for example, Ernest's (2012) search for a first philosophy of mathematics education in which he follows Levinas' argument to position the radical acceptance of the Other as its foundation.

Levinas is used to justify shifts away from the mastery of mathematical content to inherently relational concerns between people. His reception can be understood as a communicative pattern that attempts to resolve paradoxes surrounding the uniqueness of students (the "Other") in circumstances where mathematics (as a scientific discipline) and society (as the educational system's socio-political environment) assume that a universal, true and abstract body of knowledge is being transmitted. From a Laclau and Mouffean perspective, this reception seeks to displace nodal points related to the transmission of abstract mathematical knowledge in favour of nodal points of the ethical relationship to the Other.

## Calls for Hippocratic Oaths

Over the years, various people have called for Hippocratic oaths for mathematicians and educators (as reviewed in Müller et al., 2022; Rittberg, 2023). These calls are often inspired by socio-political or planetary concerns. While many of the calls originate from the applied mathematical sciences, some notable calls have also focused on education (e.g., Ardila-Mantilla, 2020; and the axioms presented in Ardila-Mantilla, 2016). Following our analytic lens, these calls can be understood as discursive interventions, whereby authors attempt to create a new nodal point (professional identity and responsibility) to re-discover and re-articulate the meaning of ethical mathematical practice in the face of today's global ethics and sustainability challenges. As a concrete example, note Freeman's (1952) call for a Hippocratic oath for statisticians, which positioned professional mathematicians and statisticians as "seeker[s] of truth in in a work environment shaped by economic and political forces" (Müller et al., 2022, p. 6):

> "Like Euclid and all the other great thinkers who have used symbols to reveal the truths of nature, I will be a seeker of the truth. Realising that numbers are only a shorthand convention for describing past events and forecasting trends, I will search for those facts expressed in numbers which show relationships and events most truly. Though surrounded by the clamour of the marketplace or of the political arena, I will not be a fraud, who selects figures to prove by chicanery a misnamed conclusion." (Freeman, 1952, p. 19)

Similarly, Sample (2019) summarises Hannah Fry's recent call as "[m]athematicians, computer engineers and scientists in related fields should take a Hippocratic oath to protect the public from powerful new technologies under development in laboratories and tech firms." Overall, we see strong socio-planetary concerns in these calls, while in the proposed codes themselves, they then translate into community concerns (What should the professional mathematician and educator do?) and mathematical concerns (What actually counts as good mathematics?). Thus, we position these calls near the "Society/Planet" vertex, but with a notable spread into the direction of both "Community" and "Mathematics".

## Positions of Professional Mathematical Societies

While the ethical and sustainable concerns of mathematical societies are slowly evolving, their codes of conduct, outputs and visions are still strongly centred on the "Mathematics" vertex, with only a slowly growing focus on community and societal concerns (cf. Buell et al., 2022; Tractenberg et al., 2024; Müller et al., 2022; Chiodo & Vyas, 2019). For example, Chiodo and Vyas (2019, p. 2) write:

> "The American Mathematical Society offers a policy statement on ethical guidelines, but this document is essentially concerned with the ethical issues that arise from mathematics as a branch of academia: aspects related to performing, publishing, and reviewing research, and those related to discrimination within the mathematical community. A similar perspective is offered by the European Mathematical Society. While the Society for Industrial and Applied Mathematics does provide a "Statement of Inclusiveness" covering discrimination, it does not appear to provide a broad ethical policy statement. The subtext here seems to be that while there are ethical issues in applied mathematics, these are imported from the disciplines that the mathematics in question is being applied to, and thus do not require a separate mention. However, their monthly publication SIAM News has, over the last decade, been publishing articles relating to ethics in mathematics which have had a gradually increasing level of depth and seriousness."

The discourses of professional societies thus reveal hegemonic articulation practices aimed at defending traditional perspectives of [M]athematics against (what is perceived as) antagonistic pressures from the community and the socio-planetary environment.

## Pragmatic Ethics in Mathematics

The pragmatic school of Ethics in Mathematics is currently primarily concerned with the prevention of harm to wider society (cf. EiMP, 2025; Müller 2024,, 2025a). As Müller (2025a) argues, this approach is largely developed *by* mathematicians *for* mathematicians in the context of higher mathematics and its education. While community concerns are certainly present (e.g., Rycroft-Smith et al., 2024; Chiodo & Bursill-Hall, 2018), the existing pragmatic approaches tend to focus less on individual student-teacher interactions or classroom interactions in general. In our analysis, this anchors much of the pragmatic scholarship between the "Mathematics" and "Society/Planet" vertices. This is not an accident, as the scholarship often explicitly tries to fill the gap between mathematical and socio-planetary concerns that would otherwise be left open. In other words, the existing pragmatic approaches represent a scientific system experiencing irritations from its environment regarding injustices produced by mathematics itself. Parts of the system have thus begun to feel a need to do something. For example, the Ethics in Mathematics Project[24] writes about its motivation that

"[n]obody can do this from outside the profession. Only mathematicians can talk to mathematicians about ethics: the discourse of philosophers is often not appropriate and specific enough for mathematicians, and simply will not address the specifically mathematical, technical ethics we face. The ethical questions of computer scientists or economists, physicists or geneticists, are pressing and important – but they are not those of the professional mathematician and they are not going to teach her how to deal with the issues that she may face in the

---

[24] Formerly known as the Cambridge University Ethics in Mathematics Project (CUEiMP).

working environment of mathematicians. So only mathematicians can do this [...]"
(EiMP, 2025)

As a related example, also consider the opening statement of the Cambridge University Ethics in Mathematics (student) Society's website, which almost entirely focuses on socio-planetary concerns:

"As such, its depth and complexity make it extremely useful. Useful because it helps us to understand the world, to our advantage. Useful because it allows us to manipulate the world, to our advantage. Useful because it empowers us to redirect the world, to our advantage. Ultimately, mathematics is useful because it gives us incredible power to change things – virtually every thing around us. But if we pause to reflect on this, we see that the utility of mathematics is derived from the way that it empowers us to understand, change, direct and manipulate the world around us, and not the other way around. It does not change the world because it is useful; it is useful because it can change the world. We see mathematics as a tool for doing good, because we can find good useful things to do with it. But none of the arguments above require us to assume that we are doing good with mathematics. It is clearly used as a way for humans to understand, change, direct and manipulate the world around us. But, just as this can be for good, it can also be for bad. Indeed, those who have the greatest ability to understand and manipulate the world hold the greatest capacity to do damage and inflict harm." (CUEiMS, n.d)

Read through the lens of Laclau and Mouffe, (some of) these existing pragmatic discourses attempt to expand the definition of mathematics by making a clearer link to the prevention of harm and re-defining the notion of good mathematics (e.g., Müller, 2024, p. 74). Finally, while critical pragmatic approaches exist (e.g., Müller, 2025b), the overarching approach has focused on re-envisioning the mathematical workflow and its education to be more inclusive of socio-planetary questions and their ethics and sustainability (e.g., Chiodo & Müller, 2025).

## Mathematics (Education) for Social Justice

The (educational) movement of Mathematics for Social Justice is positioned between community and societal concerns. Rooted in critical mathematics education, Mathematics Education for Social Justice connects classroom experiences with broader societal concerns, aiming to develop a critical consciousness in students and enabling them to understand, critique and act against social, political, and economic injustices (e.g., Buell & Shulman, 2019; Müller, 2024; Stinson et al., 2012; Karaali & Khadjavi, 2019; Gutstein & Peterson, 2005; Harrison, 2015; Shah, 2019; Bartell, 2013; Bartell et al., 2017). As a concrete example, consider Buell and Shulman's (2019) opening of a recent special issue on Mathematics for Social Justice:

> "What we teach, how we teach, and why we teach are shaped not only by institutional and accreditation requirements but also by personal philosophy. Pedagogy has three components: the curriculum, the methodology, and social education. A social justice or equity-oriented pedagogy transcends the boundaries of race, class, and gender in any classroom. In the field of mathematics, the conversation often focuses on a mathematically-rigorous curriculum (What classes should a major take? What topics are covered in Statistics?) or methodology (problem-based learning, group work, lecture, IBL, etc.). Very rarely do we discuss how to promote equity, or the role of a mathematician (or any mathematically literate person) in a democratic society. It is often argued that these issues are extracurricular ones, and do not "belong" in a mathematics classroom. However, there is a growing number of mathematicians who believe otherwise, for both pedagogical and ethical reasons." (Buell & Shulman, 2019, p. 205)

While strong mathematical concerns exist, the development of new mathematical research and knowledge is a comparatively smaller concern. Instead, a strong wish is to empower students and teachers to gain the necessary (socio-)political knowledge (e.g., Gutiérrez, 2017a), to strengthen their mathematics identity (Buell & Shulman, 2019, p. 206), and to prepare them with a (critical) general education including fostering cultural coherence, world orientation, a critical use of reason, a sense of responsibility, a strengthened sense of self, and communication and cooperation skills (cf. Heymann, 1996, p. 47), and thus, ultimately turning them into "agents of change" (Wright et al., 2024) who are able to deal with today's big social justice and sustainability problems.

From a systems-theoretic perspective, Mathematics for Social Justice can be understood as a set of communications and articulation practices which are (generally) generated by the educational system in response to strong and persistent irritations from its environment concerning socio-cultural and political injustices as well as (unethical) power dynamics (cf. discussion of Robert Moses from earlier). As is evident in Buell and Shulman's (2019) introduction to the subject, Mathematics for Social Justice effectively positions itself as a counter-hegemonic discourse that represents articulation practices which are antagonistic to the traditional (technical and abstract) views of [M]athematics, by shifting the discursive central nodal point to be social justice.

As explored earlier in the paper, educational approaches advocating for social justice and for sustainability can be related. However, here we focused on social justice to give the reader a chance to see the (subtle) differences between an approach to ethics in mathematics situated between the "Mathematics" and "Society/Planet" vertex (i.e., "Pragmatic Ethics in Mathematics") and one located between the "Community" and "Society/Planet" vertex. Their differences are further explored by Müller (2024).

## Socially-Just and Sustainable Modelling

Of the various discourses mapped, socially-just and sustainable modelling occupies the triangle's centre. This is a strategic location that enables it to function as a bridge between otherwise disparate scholarly positions. While numerous sub-discourses on modelling exist within mathematical and educational communities, for the sake of visual clarity, they are represented here as a single, sprawling discursive cloud. Its centrality stems from its inherent capacity to connect with the core concerns of all three vertices simultaneously. It speaks the language of application and rigour valued by the "Mathematics" vertex, aims to address the wish for empowerment and development of students (e.g., Heymann, 1996) central to the "Community" vertex, and engages directly with the tangible problems that define the "Society/Planet" vertex. This unique position is highlighted in survey-type discussions that showcase mathematical modelling's role in sustainable education (e.g., Karjanto, 2023; Li, 2025), and aligns with recent calls for a "critical orientation to mathematical modelling in times of disruption" (Geiger, 2024, p. 15), where the practice is not only used to solve problems but to reflexively question the assumptions and implications of the models themselves; rather than merely seeing mathematics as "embracing a positive status quote" (Makramalla, 2025, p. 539).

The connection to the "Mathematics" and "Society/Planet" vertices is particularly strong, as modelling forms the bedrock of applied mathematics (e.g., Roe et al., 2018). This makes it a familiar and acceptable practice even for those who might otherwise be sceptical of integrating ethics, challenging the NPU-assumptions from a standpoint of practical application rather than abstract critique. Discourses focused on creating balanced modelling exercises that integrate both mathematical and ethical considerations further solidify this link (e.g., Chiodo et al., 2025; Bukoski & Erbes, 2024). Simultaneously, modelling acts as a crucial nodal point for those with primarily socio-planetary concerns, providing the tools to analyse and address complex issues, and to "escape from model land" by creating models that are more socially just and contextually aware (e.g., Thompson & Smith, 2019; Thompson, 2022).

From the perspective of the "Community" vertex, socially-just and sustainable modelling is often seen as a powerful pedagogical tool for empowerment (cf. Meyer & Voigt, 2010; Meyer et al., 2025; Skovsmose 2021a, 2021b). Such modelling presents educators and students with the ability to connect the teaching and learning of mathematics with the outside world (e.g., Heymann, 1996). It shifts the focus from rote learning to active engagement, and its proponents hope that it equips students with the ability to read, understand, and critically question the models that shape their world. This critical literacy is central to learning how to both successfully "play the game" of mathematics and, more importantly, to "change the game" when its rules lead to inequitable outcomes (cf. Gutiérrez, 2009). By focusing on the development of these critical competencies, the discourse on modelling directly addresses concerns about student agency, equity, and the development of a critically conscious learning community – even though the self-efficacy necessary to raise this critical awareness and awakening (cf. Yan et al., 2024) may not always come to fruition in all students when they encounter sustainability problems, as they, among other issues, may not experience the necessary (mathematical) success (e.g., Zakariya, 2019) or the exercises are not designed carefully enough (e.g., Meyer et al., 2025).

From a discourse-theoretic perspective, the central location of socially-just and sustainable modelling allows it to function as a set of privileged articulation practices, capable of connecting differing and potentially antagonistic discourses. While the field of modelling is not without its own internal debates and tensions (Meyer & Voigt, 2010), its practices provide a common language and shared space for negotiation regarding concerns of ethics and sustainability. The articulation practices within these modelling discourses may potentially offer a pathway for scholars from different vertices of the triangle to engage in conversation, even when their foundational assumptions and ultimate goals diverge. It is this capacity to foster dialogue across divides that makes socially-just and sustainable modelling a vital and dynamic centre of gravity within the broader landscape of ethics and sustainability in mathematics and its education, and thus within the ESCT.

# Gentle and Mindful Teaching

Approaches to gentle and mindful teaching can position themselves as distinct approaches within sustainable education, usually anchored firmly within the "Community" vertex of the ESCT. Such discourses can generally be grounded on student-centred core principles, in which we see resonances with  Levinas. These can be described as the mutual respect to take one another seriously; the collaborative goal to enlighten each other; and the creative capacity to think, act, and speak in alternatives (e.g., Wilhelm, 2024, p. 144, building on Andelfinger, 1995, p. 2). As further evidenced by Wilhelm (2024), these ideas underpin an approach to mathematics education which prioritises the transformation of classroom culture and teaching practices as a prerequisite for a more sustainable world (see also Müller, 2025b). Even though it encourages Gaiatic, and thus boundary-transcending, global thinking, its primary focus remains on the immediate learning environment and the well-being and self-efficacy of learners, thus explaining its strong alignment with community-centric concerns. For example, Wilhelm writes about her conceptual approach to mindful teaching that

> "[t]he overriding goal of Mindful Teaching is to awaken the readiness of learners to use mathematics beyond the classroom to meet future challenges – for example, to participate in social decision-making in a self-determined way –, by enabling experiences of self-efficacy also in the context of sustainable development. As necessary for this purpose, Mindful Teaching values the person and the subject to a special degree – this is achieved by taking each other seriously by mathematics education and by enlightening through it." (Wilhelm, 2024, p. ii)

Analysed through a Laclau and Mouffean framework, gentle and mindful approaches to teaching represent a set of articulation practices that strategically redefine the purpose of the mathematics classroom. They work to establish student well-being, self-efficacy, and student readiness as new, fixed "nodal points," thereby challenging and seeking to displace the traditional focus on the mere transmission of abstract mathematical knowledge. While proponents of such teaching often view mathematics as a social, value-laden discipline and emphasise a globally interconnected world (as surveyed in Wilhelm, 2024), the discourse's practical and theoretical weight is applied to pedagogy and classroom dynamics. For example, one of its founders writes about gentle mathematics teaching:

> "*Gentle mathematics teaching is a culture in which Gaiatic thinking can openly engage with Cartesian-Baconian thinking. This engagement must align with the principles of 'peace', 'justice', and 'preservation of Gaia'*. [...] The presence of peace and justice in this educational concept has led me to call it 'gentle'. Gentle mathematics teaching—to put it succinctly—is achieved when it succeeds in making the described paradigmatic tensions public, addressing them, and—what is crucial—fostering affection and engagement for Gaiatic thinking." (Andelfinger, 1993, pp. 2 & 11 as cited in Wilhelm, 2024, p. 143; translated from German, italics in original)

This deliberate focus on the educational community itself is what solidifies its position near the "Community" vertex of the triangle. While making a distinction between discourses on gentle and mindful teaching and the reception of Levinas can be somewhat artificial, we positioned it slightly closer to the vertex compared to Levinas. While socio-planetary concerns may inspire such teaching (as in Wilhelm, 2024), discourses on being gentle and mindful can also exist as an ethical topic outside of classical socio-planetary sustainability concerns. Ultimately, we see the general concerns as being strongly focused on the students and other community aspects, particularly when the ethical concerns regarding being gentle and mindful are about kindness in the classroom (e.g., Baker et al., 2019).

## Ethnomathematics

Ethnomathematics connects cultural practices, societal values, and community concerns about the development and teaching of mathematics. Typically motivated by socio-political and cultural injustices, ethnomathematics examines mathematical practices globally and throughout history (D'Ambrosio, 2016; D'Ambrosio & Knijnik, 2020; Rosa & Orey, 2011; Rosa et al., 2016; Cimen, 2014; Barton, 1996; Rowlands & Carson, 2002). Here, mathematics is understood as a human, social, and culturally- and politically-situated enterprise; different mathematical practices ought to be valued for what they are. With the explicit intention of including diverse mathematical cultures, it strongly hones in on community and societal concerns, arguing for a place of ethnomathematics within mathematics. Ethnomathematically informed culturally responsive mathematics education can be used to recover lost forms of knowledge and to connect with Indigenous communities (e.g., Nicol et al, 2013).

Understood through the lens of existential sustainability, this places ethnomathematics within sustainable mathematics education and as a strong mathematical concern (Müller, 2025b). However, our analysis places ethnomathematics between the community and society vertices, albeit with a notable spread in the direction of the "Mathematics" vertex. This is because other discourses display stronger (relative) concerns regarding the development of research mathematics. This placement of ethnomathematics is perhaps the most contentious within the ESCT framework.

This placement of ethnomathematics reveals a fundamental limitation of the ESCT as a descriptive tool, as its triangular geometry reproduces the modern constitution. This constitution establishes a "great divide" between nature (seen as objective, universal, and non-human) and culture (seen as subjective, local, and human) (e.g., Goldman & Schurman, 2000; Latour, 1992, 1993, 2017, 2018). Latour (1992) says about the differences in knowledge this divide produces, that:

> "A radical is someone who claims that scientific knowledge is entirely constructed 'out of' social relations; a progressist is someone who would say that it is 'partially' constructed out of social relations but that nature somehow 'leaks in' at the end. At the other side of this tug-of-war, a reactionary is someone who would claim science becomes really scientific only when it finally sheds any trace of social construction; while a conservative would say that although science escapes from society there are still factors from society that 'leak in' and influence its development. In the middle, would be the marsh of wishy-washy scholars who add a little bit of nature to a little bit of society and shun the two extremes. This is the yardstick along which we can log most of our debates. If one goes from left to right then one has to be a social constructivist; if, on the contrary, one goes from right to left, then one has to be a closet realist. […] This tug-of-war is played in one dimension." (Latour, 1992, pp. 5-6)

The "Mathematics" vertex represents nature — a universal, a-cultural truth — while the "Community" vertex represents the culture and politics, the realm of human values and social concerns. The "Society/Planet" vertex contains aspects of both nature and culture. Ethnomathematics is a "hybrid" or "quasi-object" — a (conceptual) *monster*[25] that crosses this great divide by insisting that mathematical knowledge is a cultural practice. This is further exemplified by Rosa and Orey's (2016, pp. 11-13) six dimensions of ethnomathematics, which we summarise as follows:

- Cognitive: This dimension examines how mathematical ideas and knowledge are acquired and passed down through generations as social, cultural, and anthropological phenomena.
- Conceptual: This dimension explains how mathematical concepts originate from the practical methods and theories developed by cultural groups to meet the challenges of survival.
- Educational: This dimension seeks to humanise mathematics by incorporating values such as respect, dignity, and peace into its teaching and learning processes.
- Epistemological: This dimension investigates the diverse knowledge systems through which different cultures generate, organise, and evolve their mathematical understanding of reality.
- Historical: This dimension links the evolution of mathematics with the learner's and practitioner's own reality, examining how mathematical knowledge has been constructed and interpreted throughout history.
- Political: This dimension aims to empower students towards autonomy and citizenship by recognising and respecting the unique mathematical traditions and thinking of different cultural groups.

---

[25] In Latour's terminology monsters are not frightening or evil creatures, but actors and phenomena that defy simple categorisations through a modern perspective. The conceptual monster of ethnomathematics is thus not something "evil" or "bad," but something conceptually difficult to capture in a "modern" coordinate system. Indeed, even the "Mathematics" vertex uses the term "monster" in a similar fashion when speaking about the "monster group" in abstract algebra (consider for example, how the (mathematical) Wikipedia (2025) community talks about it: "Griess's construction showed that the monster exists.").

Consequently, the two-dimensional triangle has no fixed place for it: putting ethnomathematics into a specific location inside the triangle means ignoring its poststructuralist foundations that attempt to transcend the fixed signifiers of "nature" and "culture" altogether. The attempt to locate it on the map forces an impossible choice that reflects the logic of the modern constitution: *Is ethnomathematics really about "Mathematics" (nature), or is it really about "Community/Society" (culture)?* The difficulty in answering is not a failure of the ESCT but a success in revealing the discursive landscapes underlying political and philosophical assumptions. Ethnomathematics breaks the geometry of the triangle because it refuses the political settlement upon which the geometry of the ESCT is based. In making this visible lies the ESCT's strength. In short, ethnomathematics acts like a quantum particle when trying to be measured according to a coordinate system of classical physics; ethnomathematics breaks the modern constitution behind the triangular ontological structure of the ESCT.

In other words, the tension arises because the paper defines the "Mathematics" vertex as a formal, abstract body of knowledge, independent of human culture – what is referred to as [M]athematics rather than [m]athematics. Ethnomathematics fundamentally challenges this very idea. Its central argument is that mathematics is not a universal, pre-existing entity but rather a diverse collection of socio-culturally situated human practices, or [m]athematics. As a result, its primary concern is a critical redefinition of the "Mathematics" vertex itself. Consequently, while ethnomathematics is (from one perspective) profoundly about the nature of mathematics, its goals are not aligned with the traditional perspectives found in discourses located near the "[M]athematics" vertex. Instead, its motivations may be better understood through the other two vertices. For example, the focus on how mathematics is developed and taught within specific cultural groups is a "Community" concern. At the same time, its aim to validate non-Western knowledge systems, often in response to historical and cultural injustices, is a "Society/Planet" and "Community" concern. Therefore, positioning ethnomathematics between "Community" and "Society/Planet" is a direct result of its critical relationship with the traditionalist view that defines the ESCT's "Mathematics" vertex.

Understood through the lens of Laclau and Mouffe, the contentious location of ethnomathematics further illustrates antagonisms and hegemonic struggles over the meaning of various mathematical nodal points. Here, the nodal point of "Mathematics" (which we used to establish the first vertex of the triangle) clashes with the nodal points that ethnomathematics tries to introduce. In short, ethnomathematics could be located closer to the "Mathematics" vertex if the vertex itself were interpreted differently (see Figure 16 for a possible alternative placement). However, this location would hide tensions between classically trained mathematicians and those concerned with ethnomathematics, and thus also tensions between mathematics and approaches to its decolonisation (e.g., Borovik, 2023).

The authors also perceive this as the location effect in action. Some of the authors are research mathematicians by training, and they interpret "mathematical concerns" in a specific way that differs from many engaged with ethnomathematics. The authors place this as an invitation for conversation and further research into the ESCT. Identifying where subdiscourses self-identify, and where others place them, may shed further light on key (unspoken) tensions that the field currently experiences.

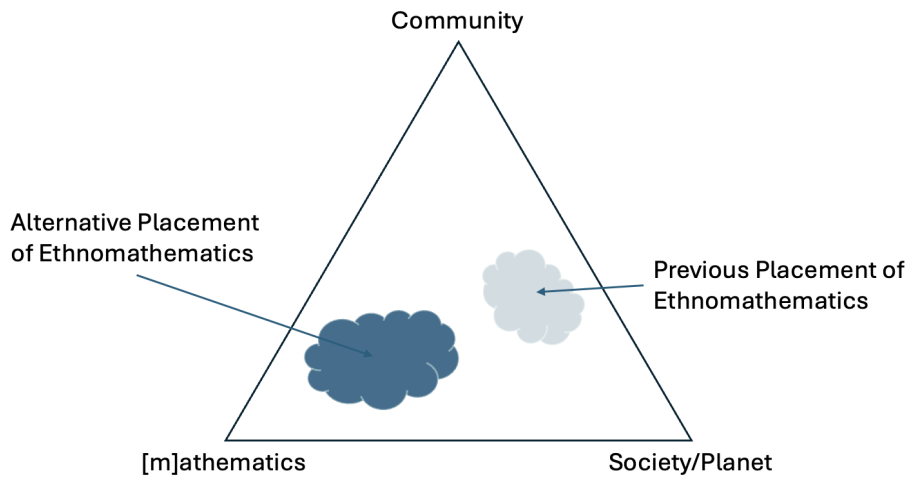

Figure 16: Alternative placement of ethnomathematics: [m]athematics vs [M]athematics

The existence of alternate placements due to the definition of the "Mathematics" vertex may not be unique for ethnomathematics, but may also hold for other discourses. If a concern is about changing (or re-defining) a vertex based on one or two of the other vertices, *our* interpretative usage of the ESCT does not necessarily capture it as a high weighting for that specific vertex. This is because the concerns, as described by Chiodo & Müller (2024), are not entirely "neutral" concerns, but they have additional qualitative dimensions. For example, the concern about "Mathematics" is a relatively "socio-politically conservative" concern to preserve the existing knowledge base and foster its growth, but the concerns about the "Community" and "Society/Planet" are more "socio-politically progressive" ("preserving what's good" and "changing what's bad").[26] Therefore, the concerns observed by Chiodo & Müller (2024) have a bias embedded within them. This bias comes to light when concerns are really concerns against, rather than concerns for, a vertex.

Our proposal of alternative placement(s) is, therefore, a significant diplomatic act, giving room for different (ethnomathematical) descriptions of mathematics (e.g., Gutiérrez's (2017b) idea of *mathematx*), and follows Barany's and Kremakova's (2023) call for a socially-just sociology of mathematics:

> "It therefore bears emphasizing, to an extent not emphasized in the foregoing analysis, that mathematics in its institutions and as a form of knowledge is today (and long has been) in many respects a tremendous locus of inequality and injustice, a resource for the empowered, and a burden for the marginalized. What

---

[26] One could also say that the formulation of the "Mathematics" vertex is more closed, while the concerns for the "Society/Planet" and "Community" vertices are more open.

we have, rather, emphasized in the foregoing is the extent to which concern for the social ramifications of mathematics has motivated and guided investigations into its social and material foundations and the potential for such investigations to revise what is taken for granted about what may and must be in the world of mathematics and from mathematics in the world." (Barany, & Kremakova, 2023, p. 18)

The map of discourses (Figure 15) reveals not just the positions of individual scholarly traditions, but also the structural dynamics of the field as a whole. We see intellectual clusters, vast distances between certain positions, and the pivotal, bridging role of discourses like socially-just and sustainable modelling. Perhaps the most telling feature of this landscape, however, is the significant discursive void between the "Mathematics" and "Community" vertices. As discussed in our analysis of different archetypal educators, this central void is not an oversight but a structural consequence of the field's dynamics. Someone situated along the "Mathematics-Community" edge has found a stable equilibrium between two attractors, namely the scientific and the educational systems, and potentially has little incentive to engage with disruptive debates on ethics and sustainability. Additionally, neither the mathematical nor the educational system appears to provide strong incentives for a nuanced middle position that argues for ethical concerns regarding both [M]athematics (with a capital M) and community. Our analysis has revealed that the nodal points of discourses near the vertices can be antagonistic (e.g., a strong concern for universal mathematical knowledge is seen as problematic by someone focused on gentle and mindful teaching; and a focus on gentle and mindful teaching could mean that students are not sufficiently quickly pushed towards the research frontier (cf. Körner, n.d.)).

Further note that while there is likely no magic educator in the middle of the triangle, there is still a discursive cloud centred there (socially-just and sustainable modelling). This is because the clouds tell us what to expect for a discourse, and the intellectual centre of mass of a discourse might not have to be near a specific educator nor near a specific paper, i.e. the "expected value" of a set of points can lie in between them, but not actually *on* any of those points. So just like the "average student" seldom exists, the "magic educator", which takes the average position of all the discourses, might not actually exist. In a way, this suggests that the perfect modelling paper also may not exist, as these also often put strong emphasis on specific concerns, e.g., the development of students' mathematical capabilities (cf. Meyer & Voigt, 2010) or the development of good mathematical models (cf. Thompson, 2022). In short, while an individual educator may find it difficult to perfectly balance all three concerns, the practice of modelling itself provides a natural meeting ground for individuals coming from all three vertices, thus placing the collective discourse on modelling in the centre.

# Conclusion

This paper introduced the Ethical and Sustainable Concerns Triangle (ESCT) to map and navigate the increasingly complex and fragmented discourses surrounding ethics and sustainability in mathematics and its education. Faced with potential communication breakdowns between different scholarly traditions, we proposed a conceptual triangle whereby its vertices represent the area's three prototypical ethics and sustainability concerns: the integrity and continuity of mathematics as a body of knowledge; the social and ethical concerns regarding the community; and broader socio-planetary concerns regarding the impact and role of mathematics on the world. The central thesis is that the fragmentation of the field can be analysed through the relational positions and dynamic tensions between these three vertices. Thus, we presented the ESCT as a meta-heuristic tool for critical reflection on the field itself.

To add a dynamic dimension to the structural map given by the ESCT, the paper continuously complemented it with a systems-theoretic interpretation. This lens reframed the triangle as an interaction between two structurally coupled systems (i.e., the scientific system of mathematics and the social communication system of education) reacting to external irritations emerging either from their shared socio-planetary environment or the other system, respectively. Discourses such as mathematics for social justice, ethnomathematics, and calls for Hippocratic oaths are thus understood as communicative responses to persistent external pressures like climate change and social injustices. The specific location of a discourse illustrates how each system processes these irritations differently. Some approaches represent a strong coupling to societal demands, leading to transformative educational goals (e.g., education for or as social justice), while others prioritise the internal, self-producing logic of their system (autopoiesis), seeking to maintain disciplinary boundaries and coherence.

A critical aspect of the ESCT is the location effect: the phenomenon whereby a scholar's or discourse's position within the triangle fundamentally shapes their perception, values, and acceptance of other scholarship. Communications originating from distant locations may be processed as noise rather than meaningful information, a communicative function rooted in differing and occasionally unstated sociomathematical and discursive norms, as well as deep-rooted philosophical assumptions about mathematics and the world. As a concrete example, the paper illustrated this effect by mapping various educator archetypes onto the ESCT, showing how their distinct positions are related to their educational and ethical priorities. This visualisation of distance between positions is not merely descriptive; it serves as a call for epistemic humility and open dialogue.

Additionally, the ESCT visualises how different systemic couplings lead to divergent priorities in school versus higher mathematics, with school-level discourses clustering more tightly around the "Community" vertex, while university-level concerns remain more oriented toward "Mathematics" and its societal applications. Similarly, the distinction between educational approaches (whether they are about, for, or as sustainability and social justice) is shown to correspond to the intensity of the systemic response to external irritations. An approach that only teaches *about* sustainability or social justice leaves the mathematical and educational systems' core logic intact, whereas an approach that strives to be education *as* social justice or sustainability represents a transformative reaction that forces a re-evaluation of the mathematical and educational systems' fundamental structures and goals.

The final mapping of key scholarly discourses provides a tangible visualisation of this entire theoretical landscape. It plots the relational positions of diverse approaches, from the community-centric reception of Levinas and the harm-prevention focus of pragmatic ethics in mathematics, to the explicitly political agenda of mathematics for social justice. This cartography makes the field's structure visible, confirming the existence of intellectual clusters, the vast distances between certain positions, and the central, bridging role of discourses like socially-just and sustainable modelling. To the authors, the contentious placement of ethnomathematics was a clear illustration of the location effect in action, revealing deep-seated struggles over the very definition of mathematics itself: Is it a formal, universal system ([M]athematics) or a collection of human, cultural practices ([m]athematics)?

Ultimately, this paper positioned the ESCT not as a final, comprehensive map, but as an analytical device and an invitation for further research and dialogue. The selection of discourses mapped is by no means exhaustive, and the framework itself is open to refinement. The authors hope that scholars and practitioners will use the ESCT as a tool for critical self-reflection, to understand their own position and its inherent perspectives and blind spots. Further use of the ESCT to locate specific discourses will hopefully continue to provide deeper insights into the interdisciplinary landscape of ethics and sustainability in mathematics and its education.

**AI Usage Statement:** The authors utilised artificial intelligence to support the preparation of this manuscript. Specifically, Gemini Deep Research and Elicit were used to assist with the literature search. For improving grammar, clarity, and language, the authors employed Grammarly, Gemini, ChatGPT, and Claude. The authors reviewed and edited all AI-assisted content and take full responsibility for the final text.

# Appendix A: "Pragmatic Ethics in Mathematics"

**Summary:** To illustrate the discourse-analytic mapping process, this appendix shows how Chiodo & Müller analysed their own pragmatic approach using the methodology displayed in figures 2 - 4.

| Educational Context | Mathematics (Concerns) | Community (Concerns) | Society / Planet (Concerns) |
|---|---|---|---|
| **Authorial Stance** | strong | weak | strong |
| **School of Thought** | Strong (pragmatic, classical) | medium | medium to strong |
| **Expressed and Hidden Concerns** | strong | weak | strong |

| Identities, Hegemonies and Power Relations | Mathematics (Concerns) | Community (Concerns) | Society / Planet (Concerns) |
|---|---|---|---|
| **Levels of Ethical Engagement** | medium | medium | strong |
| **Core Assumptions** | strong | weak | strong |
| **Philosophical Foundation** | strong | weak | medium |
| **Educational Position** | medium | medium | medium |

| Identification of Mathematical Nodal Points | Pillars |
|---|---|
| **Strong engagement** | 1, 3, 4, 5, 9, 10 |
| **Medium engagement** | 6, 7, 8 |
| **Weak engagement** | 2 |

Table 5: Discursive concerns for "Pragmatic EiM"

# Appendix B: Index of the Core Ideas of the Paper

| Component | Concept | Function |
|---|---|---|
| The three vertices (mathematics, community, society/planet) | Each vertex represents a core concern: the ethics and sustainability concerns of ensuring that research mathematics remains a stable and continuing body of knowledge; the ethical, sustainability, and social questions surrounding our mathematical, and other, related communities; and questions regarding mathematics within and its impact on the wider society and our planet. | The vertices act as orienting poles. The proximity to a vertex indicates the degree of focus that is put on a specific concern. All locations in the triangle are convex combinations of concerns. The more weight a concern is given in a discourse, the closer it is situated to the corresponding vertex. This can be understood to be a locational grammar of discourses on ethics and sustainability. |
| The location effect | A scholar's or discourse's position within the triangle influences their perception, interpretation, and acceptance of other scholarship. Discourses that are close to each other, can usually effectively engage with each other. Discourses far away from each other may face strong tensions, and could be interpreted as incompatible by some. | The location effect explains some of the field's fragmentations. Communications (e.g., knowledge claims) far away from one's own position may be perceived as noise. The location effect highlights the need for epistemic humility and open dialogue. |
| Systems perspective | By changing to a systems theoretic perspective, the static (ESCT) triangle can be interpreted as a dynamic interaction between two systems and their environment. The vertex "Mathematics" represents a subsystem of science, and the "Community" vertex represents the social communication system of education; both systems exist in a shared socio-planetary environment ("Society/Planet" vertex) sending irritations to the systems, and they also irritate each other. | The systems theoretic interpretation explains the emergence of different discourses as reactions to irritations (e.g., evidence of climate change). For example, injustices in the socio-planetary environment (or within the scientific subfield of mathematics) can lead to reactions in the educational system, which begins to respond by advocating for curricular changes, new pedagogical strategies, etc. |
| Educator Archetypes | The five educator archetypes described by Ernest (1991): Old humanist mathematicians, progressive educators, public educators, technological pragmatists, and industrial trainers. | The different educators represent different couplings between the education system and its environment (e.g., the scientific system or the wider economic, socio-political environment) |
| NPU (Neutrality, Purity, Universality) Assumptions | The NPU Assumptions represent various assumptions surrounding the neutrality, purity and universality of mathematics. | The NPU assumptions can be represented as a gradient across the triangle: increasing distance from the "Mathematics" vertex is correlated to decreasing NPU. |
| Education about, for or as sustainability and social justice | These three approaches represent different orders of learning (first-, second- and third-order learning); they are related to different levels of engagement and to perspectives about how mathematics, education, and society have to change to deal with, for example, climate change. | They can be understood as a proxy to signal how strongly environmental pressures have acted on mathematics and the education system: "about" leaves structures largely intact; while "for" and "as" ask for increasingly holistic changes. |
| Discourses and their discursive clouds | This is the visual map of locations and spreads of different scholarly discourses. The discourse clouds are a metaphorical visualisation of their spread of concerns. | The clouds represent relational positions (and thus choices and tensions) that structure ethics and sustainability in mathematics and its education. |

Table 6: Summary of core ideas

# Appendix C: Guidance for Self-Positioning

While the earlier tables 2 - 4 explained how the authors located different discourses inside the triangle, the following table can help the reader to position themselves inside the triangle. The questions are not exhaustive, but provide a starting point for what to ask.

| Vertex | Example Questions | Orientation | Systems-Theoretic Aspect |
|--------|-------------------|-------------|--------------------------|
| **Mathematics** | Do I see mathematics as neutral (i.e., value-free), pure (i.e., abstract), and universal? Is preserving mathematical knowledge, rigour or elegance one of my central concerns? | Potential focus on mathematical purity, classical philosophies of mathematics (Platonism, Formalism, Logicism) | Strong emphasis on autopoiesis of the mathematical system. Potentially a low sensitivity to irritations from society (or from education), much of it may be perceived as noise. |
| **Community** | Do I focus strongly on the well-being and personal development of learners? Am I concerned with issues of diversity, equity, inclusion, identity, etc.? Do I view mathematics education as situated in a socio-cultural context? | Potential emphasis on pedagogy, issues of care, interpersonal responsibility, teacher-student interactions, etc. | Strong emphasis on the autopoiesis of the education system. Potentially strongly coupled with the learners' inner psychic systems, that may trigger new concerns. |
| **Society / Planet** | Am I particularly concerned about using mathematics to address societal or planetary concerns? Do I believe that education should serve a higher, external societal goal (e.g., developing a sustainable society at large)? | Potential for using critical mathematics education, education for/as sustainability/social justice; potential to strongly advocate for student agency in society that uses and questions mathematics and how its deployed | Strong sensitivity (i.e., awareness) to receiving irritations from the larger socio-planetary environment |
| **Between Vertices** | Do I feel torn between, or have affinity/sympathy for, different ethical and sustainable concerns? Am I constantly trying to translate between different groups and demands? | Potential focus on mixed areas (e.g., sustainable and socially just modelling; mathematics for social justice, etc.) | Individual who is strongly coupled with multiple systems and/or the environment; potentially experiences a lot of inner or external conflicts regarding their work/research/pedagogy |
| **Dissonance and Shifts** | Are there areas that I find difficult to understand? Are there areas that I understand but cannot accept? Have my views on mathematics, community or society recently changed? | This suggests that an individual is currently dynamically moving within the triangle (e.g., someone who has just shifted their research, moved to a different teaching position, experienced or learned about something new, etc.) | Dissonances and shifts are indicators for systemic reconfiguration. They likely represent a growing awareness according to Rycroft-Smith et al.'s (2024) levels of ethical engagement. External irritations may be currently perceived in a more welcoming fashion. |

Table 7: Guidance for self-positioning